\documentclass[11pt,twoside]{article}
\usepackage{atmp}
\copyrightnotice{2003}{7}{87}{120}	

\newcommand{\bea}{\begin{eqnarray}}
\newcommand{\ena}{\end{eqnarray}}
\newcommand{\be}{\begin{eqnarray*}}
\newcommand{\en}{\end{eqnarray*}}

\newcommand{\lb}[1]{\label{#1}}
\newcommand{\noi}{\noindent}

\newcommand{\dstyle}{\displaystyle}

\newcommand{\N}{{\mathbb N}}
\newcommand{\R}{{\mathbb R}}
\newcommand{\Z}{{\mathbb Z}}
\newcommand{\C}{{\mathbb C}}
\newcommand{\calC}{{\mathcal C}}
\newcommand{\F}{{\mathcal F}}
\newcommand{\Fp}{{\mathcal F}^\phi}
\newcommand{\Fg}{{\mathcal F}^{\beta\gamma}}
\renewcommand{\H}{{\mathcal H}}

\newcommand{\Phibar}{\bar{\Phi}}
\newcommand{\Jbar}{\bar{J}}
\newcommand{\Lambdabar}{\bar{\Lambda}}

\newcommand{\Mhat}{\widehat{M}}
\newcommand{\Lhat}{\widehat{L}}
\newcommand{\Jhat}{\widehat{J}}
\newcommand{\Phat}{\widehat{P}}
\newcommand{\Qhat}{\widehat{Q}}
\newcommand{\Rhat}{\widehat{R}}
\newcommand{\What}{\widehat{W}}

\newcommand{\hhat}{\hat{\mathfrak{h}}}
\newcommand{\Deltahat}{\widehat{\Delta}}
\newcommand{\deltahat}{\hat{\delta}}
\newcommand{\alphahat}{\hat{\alpha}}
\newcommand{\lambdahat}{\hat{\lambda}}
\newcommand{\muhat}{\hat{\mu}}
\newcommand{\rhohat}{\hat{\rho}}

\newcommand{\h}{\mathfrak{h}}

\newcommand{\tr}{{\rm tr}}
\newcommand{\ch}{{\rm ch}}
\newcommand{\Hom}{{\rm Hom}}

\newcommand{\bra}[1]{\langle #1 |}        
\newcommand{\ket}[1]{{| #1 \rangle}}      

\newcommand{\sltwohat}{\widehat{\mathfrak{sl}}_2}
\newcommand{\slNhat}{\widehat{\mathfrak{sl}}_N}

\newcommand{\slN}{\mathfrak{sl}_N}

\theoremstyle{plain}
\newtheorem{thm}{Theorem}[section]
\newtheorem{cor}[thm]{Corollary}
\newtheorem{prop}[thm]{Proposition}
\newtheorem{lem}[thm]{Lemma}
\newtheorem{define}[thm]{Definition}

\theoremstyle{remark}
\newtheorem*{rem}{Remark}

\begin{document}

 \title{Affine Analogue of Jack's Polynomials for $\sltwohat$}
 
 \url{math.QA/0210236}		
 
 \author{Yuji Hara}		
 \address{Graduate School of Mathematical Sciences\\ University of Tokyo\\Tokyo 153-8914, Japan}
 \addressemail{snowy@ms.u-tokyo.ac.jp} 	
 
 \markboth{\it Affine analogue of Jack's \ldots}{\it Yuji Hara}

 \begin{abstract}
Affine analogue of Jack's polynomials introduced by Etingof and Kirillov Jr. is studied for the case of $\sltwohat$.
Using the Wakimoto representation, we give an integral formula of elliptic Selberg type for the affine Jack's polynomials.
From this integral formula, the action of $SL_2(\Z)$ on the space of the affine Jack's polynomials is computed.
For simple cases, we write down the affine Jack's polynomials in terms of some modular and elliptic functions.
 \end{abstract}

\cutpage

\section{Introduction}
In recent years,  beautiful connections between the theory of symmetric functions and other branches of mathematics and physics have been studied. Those branches include $q$-deformed spherical functions \cite{M1}\cite{N}, quantum-many body problems e.g the Calogero-Moser-Sutherland models (see \cite{CMS} and references therein), representation theory of affine Hecke algebras \cite{C} and representation theory of Lie algebras and quantum groups \cite{EK}. Among many symmetric functions, Macdonald's polynomials and Jack's polynomials \cite{M} have attracted so much attention which revealed the rich mathematical structure behind them.

In their theory connecting those polynomials to the representation theory of Lie algebras and quantum groups, P. Etingof and A. Kirillov Jr. introduced an affine analogue of Jack's symmetric functions. These functions are called affine Jack's polynomials even though they are not polynomials but a generalization of characters of affine Lie algebras.
They first introduced differential operators associated to affine root systems of type $A^{(1)}_{N-1}$.  They are an elliptic generalization of the Calogero-Sutherland type hamiltonians except that they contain a term of differentiation w.r.t the elliptic nome.
The affine analogue of Jack's polynomial $\Jhat_{\lambdahat}$ is labeled by the dominant integral weight $\lambdahat$ and defined by two conditions: (1) it is an eigenfunction of the affine Calogero-Sutherland type operator with a certain eigenvalue, (2) its transition matrix w.r.t affine analogue of monomial symmetric functions is upper triangular.
The very important property of this affine Jack's polynomials is that they can be given as traces of vertex operators over highest weight modules:\,Theorem \ref{thm:J-tr}.
Two more of three important aspects of the affine Jack's polynomials can be seen from this property:
\begin{itemize}
\item[(1)]
Affine analogue of Jack's symmetric functions.
\item[(2)]
One point functions of Conformal Field Theory on a torus. The differential operator is nothing but the Knizhnik-Zamolodchikov-Bernard heat operator.
\item[(3)]
Affine Jack's polynomials can be considered as a generalization of characters of affine Lie algebras. They form a basis of a space invariant under the action of $SL_2(\Z)$
\end{itemize}

The aim of this paper is to give explicit expressions for the affine Jack's polynomials itself and their modular transformation property in the case of $\sltwohat$. In Section \ref{sec:intgrl}, we give an integral formula of the affine Jack's polynomials. This formula is calculated by means of the Wakimoto representation of $\sltwohat$. In Section \ref{sec:modular}, using this integral formula, we computed the action of  $SL_2(\Z)$ on the space of the affine Jack's polynomials. And this projective representation of $SL_2(\Z)$ allows us to write down the affine Jack's polynomials in terms of some elliptic and modular functions in Section \ref{sec:example}. Together with the integral formula given in Section \ref{sec:intgrl}, this gives some formulas for elliptic Selberg type integrals.

In \cite{FSV}, similar results are given for some solutions of the Knizhnik-Zamolodchikov-Bernard heat equation. The main differences are that their approach is geometrical and not considering the affine Jack's polynomials specially. Also, during the preparation of this paper we found \cite{FSV2} which gives a result similar to our result presented in Section \ref{sec:example}.

\section{Affine Analogue of Jack's Polynomials}
We fix our notations for $\slNhat$ \cite{Ka} and review the affine Jack's polynomials. See \cite{EK} for the detail. In this section, we consider $N=2,3,4 \cdots$.

\subsection{Affine Lie algebras}
Let us fix the notation for $\slN$ first.
Let $\h$ be the Cartan subalgebra and $\h^*$ be its dual.
The standard invariant bilinear form on $\slN$ is denoted by $(\cdot,\cdot)$ and normalized so that the induced form on $\h^*$ satisfies $(\alpha,\alpha)=2$ for a root $\alpha$.
Let $\alpha_1, \ldots,\alpha_{N-1}$ be simple roots and the longest root be $\theta$.
Let $R$ be the root system. We fix a polarization $R=R^+\sqcup R^-$ where $R^{\pm}$ are the subsets of positive and negative roots respectively.
We denote the fundamental weights by $\Lambdabar_i\in\h^*\,(i=1,\cdots,N-1)$.
Weyl vector is given by $\rho = \sum_{i=1}^{N-1} \Lambdabar_i$.
Let $P = \{ \lambda\in\h^* | (\lambda,\alpha_i) \in\Z \} = \bigoplus_i \Z\Lambda_i$ be the weight lattice and $P^+ = \{ \lambda\in\h^* | (\lambda,\alpha_i) \in\Z_+ \} = \bigoplus_i \Z_+\Lambda_i$ be the cone of dominant weights.

\noi
The affine Lie algebra $\slNhat$ can be realized as an extension of a loop algebra:
\begin{gather*}
\slNhat = \slN \otimes \C[t,t^{-1}] \oplus \C c \oplus\C d,
\\
[x\otimes t^n,y\otimes t^m] = [x,y] \otimes t^{m+n}+n \delta_{m,-n}(x,y)c,
\\
c  \text{ is central},\quad
[d, x\otimes t^n] = nx \otimes t^n.
\end{gather*}
We use the notation $x[n] = x \otimes t^n$. Sometimes we use a smaller algebra $\slNhat' = \slN \otimes \C[t,t^{-1}] \oplus \C c$.

Cartan subalgebra and its dual are given by
$\hhat = \h \oplus \C c \oplus \C d, \hhat^* = \h^* \oplus \C \delta\oplus \C\Lambda_0.$
The canonical pairing $\langle\cdot,\cdot\rangle$ between $\hhat$ and $\hhat^*$ is given by
$\langle\Lambda_0, \h\oplus\C d\rangle=\langle\delta, \h\oplus\C c\rangle = 0, \langle\delta,d\rangle = 1,\langle\Lambda_0, c\rangle = 1.$
Bilinear non-degenerate symmetric form on $\hhat^*$ is given by
$(\Lambda_0, \delta)=1, (\Lambda_0,\h^*)=(\delta,\h^*)=(\Lambda_0,\Lambda_0)=(\delta,\delta)=0.$
This pairing and the bilinear form coincide with those of $\slN$ on $\h$ and $\h^*$.
Let
$\Rhat = \{\alphahat = \alpha+n\delta\,|\,\alpha\in R, n\in\Z\text{ or }\alpha=0, n\in\Z\setminus\{0\}\}$
be the affine root system.
We fix the polarization by taking the set of positive roots as
$\Rhat^+ = \{\alphahat = \alpha + n\delta \in \Rhat \,|\, n>0 \text{ or } n=0, \alpha\in R^+\}$.
Simple roots are given by $\alpha_0 = -\theta+\delta, \alpha_1, \ldots,\alpha_{N-1}.$
Let
$\hat{Q} = \bigoplus \Z\alpha_i$
be the affine root lattice and
$\hat{Q}^+ = \bigoplus \Z_ + \alpha_i$
be the positive part.
Let $\Lambda_i\,(i=0,\cdots ,N-1)$ be the fundamental weights where
$\Lambda_i=\Lambdabar_i +\Lambda_0$ for $i=1,\cdots ,N-1.$
The affine weight lattice is given by
$\Phat = P \oplus \Z\delta \oplus \Z\Lambda_0 \subset \hhat^*$
and the cone of dominant weights is
$\Phat^+ = \{ \lambdahat\in\hhat^* \,|\, \langle\lambdahat,\alpha^\vee_i\rangle \in \Z_+, i=0,\ldots,r \}.$
For $K\in\Z$, we define
$\Phat^+_K = \{ \lambda + n\delta + K\Lambda_0 \,|\, \lambda\in P^+ \}$
where
$P^+_K = \{\lambda\in P^+ \,|\, (\lambda,\theta)\leq K\}.$
Affine Weyl vector is defined as $\rhohat = \sum_{0=1}^{N-1} \Lambda_i$.
We set an order in $\Phat$ by
$\lambdahat \leq \muhat \iff \muhat - \lambdahat \in \Qhat^+. $

Let $\C[\Phat] = \bigoplus_{K\in\Z} \C[\Phat_K]$ be the group algebra of the weight lattice.
We consider a completion:
\begin{gather*}
\overline {\C[\Phat]} = \bigoplus_K \overline{\C[\Phat_K]},
\\
\overline{\C[\Phat_K]} =
  \left\{ \sum_{ \lambdahat\in\Phat_K } a_{\lambdahat} e^{\lambdahat} \biggm|
  {}^\exists\! N \mbox{ s.t. } a_{\lambdahat} =0 \mbox{ for } (\lambdahat, \rhohat)\geq N \right\}.
\end{gather*}
This completion is chosen so to include the characters of the Verma modules (and more generally, modules from the category $\mathcal{O}$) over $\slNhat$. Let $A$ be the subalgebra of $\overline {\C[\Phat]}$ defined by
\begin{align*}
A = \bigcap_{w\in\What} w\biggl(\overline{\C[\Phat]}\biggr),
\quad
A_K = \bigcap\limits_{w\in\What} w\biggl(\overline{\C[\Phat_K]}\biggr).
\end{align*}
Let $\What$ be the affine Weyl group of $\slNhat$. One of the main objects of our study will be the algebra of $\What$-invariants $A^{\What}$.
We introduce a formal variable $p=e^{-\delta}$. Then every element of $A$ can be written as a formal Laurent series in $p$ with coefficients from $\C[P]$.
The following properties characterize $A$ as the affine analogue of $\C[P]$.
\begin{lem}
For any $\lambdahat\in \Phat_K, K\geq 0$ the orbitsum
\begin{align*}
m_{\lambdahat} = \sum_{\muhat\in\What\lambdahat} e^{\muhat}
\end{align*}
belongs to $A^{\What}_K$.
\end{lem}

\begin{lem}
If $M$ is a module from the category $\mathcal{O}$ then $\ch M\in A$ iff $M$ is integrable, in which case $\ch M\in A^{\What}$. In particular, the irreducible quotient $\ch L_{\lambdahat} \in A^{\What}$ iff  $\lambdahat\in \Phat^+$.
\end{lem}

\begin{lem}
$A^{\What}_K=0$ for $K<0$, and $A^{\What}_0 = \left\{\sum_{n\leq n_0} a_n e^{n\delta} | a_n\in\C\right\}.$
\end{lem}

\begin{thm}
For every $K\in\Z_+$, the orbitsums $m_{\lambda + K\Lambda_0}, \lambda\in P^+_K$ form a basis of $A^{\What}_K$ over the field $\C((p))$.
\end{thm}

We also need $\hat{\mathcal{R}}$ which is an extension of $A$.
Let us consider the algebra $\C[\Phat](1-e^{\alphahat})^{-1}$, obtained by adjoining to $\C[\Phat]$ the inverse of the form $(1-e^{\alphahat})$:
\begin{align*}
\C[\Phat](1-e^{\alphahat})^{-1} =
\left\{ \frac{f}{g} \biggm| f,g\in\C[\Phat], \ g=\prod_{\alphahat\in I} (1-e^{\alphahat}) \mbox{ for } I\subset \Rhat \right\}.
\end{align*}
We have a morphism for every $w\in\What$:
\begin{align*}
\tau_w\ :\ & \C[\Phat](1-e^{\alphahat})^{-1} \to w\left(\overline {\C[\Phat]}\right),
\\
& (1-e^{\alphahat})^{-1} \mapsto \sum_{n=0}^\infty e^{-n\alphahat} \quad\mbox{ for } \alphahat\in w\Rhat^+.
\end{align*}
Then $\hat{\mathcal{R}}$ is given by
\begin{align*}
\hat{\mathcal{R}} =
\biggl\{
  \sum a_n \bigm|  a_n\in\C[\Phat](1-e^{\alphahat})^{-1}, \;
  \sum \tau_w(a_n) \in w\biggl(\overline{\C[\Phat]}\biggr) \mbox{ for } {}^\forall\!w\in \What
\biggr\}.
\end{align*}

\subsection{Affine Analogue of Jack's Polynomials}
\label{subsec:affine Jack}
Throughout this paper we consider parameters $K$ and $k$ (introduced below) as $K\in\Z,\;k\in\N$.

Let us introduce the following differential operators, which we consider formally as derivations of the algebra $\C[\Phat]$:
\begin{gather*}
\Deltahat e^{\lambdahat} = (\lambdahat, \lambdahat) e^{\lambdahat},
\\
\partial_{\alphahat} e^{\lambdahat} = (\alphahat, \lambdahat) e^{\lambdahat}, \quad \alphahat \in \hhat^*.
\end{gather*}
If we use the notation $p=e^{-\delta}$ and write elements of $\C[\Phat_K]$ as functions of $p$ with coefficients from $\C[P]$: $\C[\Phat_K] = e^{K\Lambda_0} \C[p, p^{-1}][P]$ then
\begin{gather*}
\Deltahat = -2K p\frac{\partial}{\partial p} +\Delta_{\h},
\\
\partial_{\alpha + n\delta} = \partial_\alpha + nK
\end{gather*}
where $\Delta_\h$ is the Laplace operator on $\h$:
\begin{align*}
\Delta_{\h} e^{\lambda} = (\lambda, \lambda) e^{\lambda},\;\lambda\in\h^*.
\end{align*}
Let
$\deltahat = e^{\rhohat} \prod_{\alphahat \in \Rhat^+} (1-e^{-\alphahat})$
be the affine Weyl denominator.
\begin{define}
The Calogero-Sutherland operator $\Lhat_k$ for the affine root system $\Rhat$ is the differential operator which acts on $\hat{\mathcal{R}}_K$:
\begin{align}
\Lhat_k = \Delta - 2Kp\frac{\partial}{\partial p}
                - k(k-1)
				\sum_{\begin{subarray}{c}
				\alpha \in R^+ \\ n\in\Z
				\end{subarray}}
				\frac{p^n e^\alpha}{(1-p^n e^\alpha)^2} (\alpha,\alpha).
\label{eqn:L}				
\end{align}
We mainly use the operator $\Mhat_k$ conjugated by the affine Weyl denominator:
\begin{align}
\Mhat_k &= \deltahat^{-k}\circ (\Lhat_k - k^2 (\rhohat,\rhohat) ) \circ\deltahat^{k}
\\
&= \Deltahat - 2k \sum_{\alphahat \in \Rhat^+} \frac{1}{1-e^{\alphahat}} \partial_{\alphahat} + 2k \partial_{\rhohat}.
\end{align}
\end{define}
\begin{define}
\label{def:affineJack}
Affine Jack's polynomial $\Jhat_{\lambdahat},\; \lambdahat\in\Phat^+$ is the element of $A^{\What}$ defined by the following conditions:
\begin{itemize}
\item[(1)]
${\displaystyle
\Jhat_{\lambdahat} = m_{\lambdahat} + \sum_{\muhat < \lambdahat} c_{\lambdahat,\muhat} m_{\muhat},
}$
\item[(2)]
${\displaystyle
\Mhat_k \Jhat_{\lambdahat} = (\lambdahat, \lambdahat + 2k \rhohat) \Jhat_{\lambdahat}.
}$
\end{itemize}
\end{define}

\begin{rem}
We use the terminology \lq polynomial' following \cite{EK} even though it is not a polynomial but a theta function.
\end{rem}

We denote by $L_{\lambdahat}$ an irreducible highest weight representation of the highest weight $\lambdahat\in\Phat$.
Let $V$ be the finite dimensional irreducible representation of $\slN$. The evaluation module based on $V$ is defined as
\begin{gather*}
V(\zeta)=V\otimes\C[\zeta,\zeta^{-1}],
\\
\pi_{V(\zeta)}(a[n]) = \zeta^n \pi_{V}(a),
\\
\pi_{V(\zeta)}(c) = 0,\quad \pi_{V(\zeta)}(d) = \zeta\frac{d}{d\zeta}.
\end{gather*}

Let $S^{N(k-1)}\C^N$ be a symmetric tensor of the vector representation of $\slN$.
Let $\muhat\in\Phat^+$. A non-zero intertwiner of $\sltwohat'$
\begin{align}
\Phi\,:\, L_{\muhat} \mapsto L_{\muhat}\otimes \left( S^{N(k-1)}\C^N \right) (\zeta)
\label{eqn:Phi}
\end{align}
exists iff $\muhat = \lambdahat + (k-1)\rhohat,\; \lambdahat\in\Phat^+$;  if it exists, it is unique up to a constant. We will denote such an intertwiner by $\Phi_{\lambdahat}$.
\begin{define}
\begin{align}
\varphi_{\lambdahat} = \sum_{\muhat} e^{\muhat} \,\tr|_{L_{\lambdahat}[\muhat]}(\Phi_{\lambdahat}).
\end{align}
\end{define}
It takes values in the weight zero subspace $\left( S^{N(k-1)}\C^N \right) \![0]$. This subspace is one-dimensional and can be identified with $\C$. This identification and the normalization of $\Phi_{\lambdahat}$  are so chosen that $\varphi_{\lambdahat} = e^{\lambdahat + (k-1)\rhohat} + \cdots$. Hereafter we basically consider $\varphi_{\lambdahat}$ as a scalar function in this way.
\begin{thm}
\label{thm:J-tr}
\begin{itemize}\item[]
\item[(1)]
${\displaystyle
\varphi_0 = \deltahat^k,
}$
\item[(2)]
${\displaystyle
\Jhat_{\lambdahat} = \frac{\varphi_{\lambdahat}}{\varphi_0}.
}$
\end{itemize}
\end{thm}
\noi
From this formula it is easy to see that $\Jhat_{\lambdahat + n\delta} = p^{n} \Jhat_{\lambdahat}$.

\vspace{5mm}
Let us consider the following domain
\begin{align*}
Y = \h \times \C \times \H
\end{align*}
where $\H$ is the upper half-plane: $\H=\{\tau\in\C|\text{ Im }\tau >0\}$.
Then every element $e^{\lambdahat} \in \C[\Phat]$ can be considered as a function on $Y$ as follows: if
$\lambdahat = \lambda + a\delta + K\Lambda_0$ then put
\begin{align*}
e^{\lambdahat} (h,u,\tau) = e^{2\pi i[\langle\lambda,h\rangle+Ku-a\tau]}.
\end{align*}
Note that this in particular implies that $p=e^{-\delta}$ is given by $p=e^{2\pi i\tau}.$

It is easy to see that if $f\in \C[\Phat]$ then $f(h, u, \tau+1)=f(h,u, \tau)$, so we can as well consider $f$ as a function of $h, u, p$, where $p\in \C$ is such that $0<|p|<1$. Note that $f \in \C[\Phat_K]\iff f(h, u,\tau)=e^{2\pi i K u}f(h,0, \tau)$; in this case we say that $f$ is a function of level $K$.
\begin{thm}
For every $\lambdahat \in \Phat_K$, $\Jhat_{\lambdahat}$ and $\varphi_{\lambdahat}$ converge on $Y$.
\end{thm}

Now let us recall some facts about the modular group and its action. The modular group $\Gamma = SL_2(\Z)$ is generated by the elements
\begin{align*}
S =
\begin{pmatrix} 0 & -1 \\ 1 & 0 \end{pmatrix},
\quad
T =
\begin{pmatrix} 1 & 1   \\ 0 & 1 \end{pmatrix}
\end{align*}
satisfying the defining relations $(ST)^3=S^2, S^2T=TS^2, S^4=1$.
This group acts in a natural way on $Y$ as follows:
\begin{align*}
\begin{pmatrix} a & b \\ c & d \end{pmatrix} (h,u,\tau)
= \left( \frac{h}{c\tau+d} , u-\frac{c(h,h)}{2(c\tau+d)} , \frac{a\tau+b}{c\tau+d} \right).
\end{align*}
In particular,
\begin{align*}
T(h,u,\tau) &= (h, u,\tau+1),
\\
S(h,u,\tau) &= \left(\frac h \tau , u-\frac{(h,h)}{2\tau} ,-\frac{1}{\tau}\right).
\end{align*}
Also, for any $j\in \C$ we will  define a right action of $\Gamma$ on functions on $Y$ as follows: if
$\alpha = \begin{pmatrix} a & b \\ c & d \end{pmatrix}$ then let
\begin{align}
(f[\alpha]_j) (h,u,\tau) = (c\tau+d)^{-j} f(\alpha(h,u,\tau)).
\end{align}
This is a projective action, which is related to the ambiguity in the choice of $(c\tau+d)^{-j}$ for non-integer $j$.
\begin{define}
Let $K\in\Z_+,\,\lambda\in P^+_K$. Normalized affine Jack's polynomial $J_{\lambda,K}$ is defined by
\begin{gather}
J_{\lambda,K} = \Jhat_{\lambda+\alpha\delta+K\Lambda_0},
\\[2mm]
\alpha = \frac{k(\rho,\rho)}{2h^\vee} - \frac{(\lambda+k\rho,\lambda+k\rho)}{2(K+kh^\vee)}
\end{gather}
where $h^\vee$ is the dual Coxeter number of \ $\slNhat$.
\end{define}
\noi
{\it Remark}
\begin{gather}
J_{\lambda,K} =
\frac{1}{ \deltahat^{\prime\, k-1} }
\tr_{L_{ \lambda + K\Lambda_0 + (k-1)\rhohat }} \left( \Phi_{\lambda+\alpha\delta+K\Lambda_0}\, p^{ L_0 - \frac{c_V}{24} } e^{2\pi i h} \right),
\\[2mm]
\deltahat^{\prime\, k-1}
= \tr_{L_{ (k-1)\rhohat }} \left( \Phi_{\lambda+\alpha\delta+K\Lambda_0}\, p^{ L_0 - \frac{c_V}{24} } e^{2\pi i h} \right),
\quad
\deltahat^\prime = ( p^{ \frac{(\rho,\rho)}{2h^\vee} } \deltahat )^{k-1}
\end{gather}
where $c_V$ is the central charge for the action of the Virasoro algebra on $L_{ \lambda + K\Lambda_0 + (k-1)\rhohat }$:\, $c_V = (K+(k-1) h^\vee )\, {\rm dim}\,\slN  / (K+ kh^\vee)$.
\begin{thm}
Let $K\in \Z_{\geq 0}$. Denote
\begin{align*}
V_K = \bigoplus_{\lambda\in P^+_K} \C\,J_{\lambda,K}.
\end{align*}
Consider elements of $V_K$ as functions on $Y$. Then $V_K$ is preserved by the action of $\Gamma$ with weight $j = -\frac{K(k-1)r}{2(K+kh^\vee)}$. In particular, this means that $V_K$ is naturally endowed with a structure of a projective representation of $\Gamma$.
\end{thm}

\section{Integral formula via Wakimoto representation}
\label{sec:intgrl}
Hereafter we restrict our attention to $\sltwohat$.

\subsection{Wakimoto representation}
Let us introduce the bosonic ghost $\beta_n, \gamma_n\;(n\in\Z)$, the boson $\phi_n\;(n\in\Z),Q$ and the degree operator $L_0$ by
\begin{gather*}
[\beta_n,\gamma_m]=\delta_{n+m,0},
\\
[\phi_n,\phi_m]=\frac{2n}{\kappa}\delta_{n+m,0},\quad [\phi_0,Q]=\frac{2}{\kappa},
\\
\quad [L_0, b_n] = n b_n \quad (b = \beta, \gamma, \phi).
\end{gather*}
Other combinations are commutative.
Currents are defined as
\begin{gather*}
\beta(\zeta) = \sum_{n\in\Z}\beta_n \zeta^{-n-1},\quad \gamma(\zeta)=\sum_{n\in\Z}\gamma_n \zeta^{-n},
\\
\phi(\zeta) = Q + \phi_0\ln \zeta - \sum_{n\neq 0}\frac{\phi_n}{n}\zeta^{-n},
\quad
\partial\phi(\zeta) = \sum_{n\in\Z}\phi_n\zeta^{-n-1}.
\end{gather*}
We define the Fock spaces $\F_J, \Fg, \Fp_J$ by
\begin{align*}
&\F_J=\Fg\otimes\Fp_J,
\\
&\Fg=\C[\beta_{-1},\cdots,\gamma_{0},\cdots]\ket{vac},
\\
&\beta_n\ket{vac}=0\quad(n\geq0),\quad\gamma_n\ket{vac}=0\quad(n>0),
\\
&\Fp_J=\C[\phi_{-1}\cdots]\ket{2J+1},
\\
&\phi_0\ket{2J+1}=\frac{2J}{\kappa}\ket{2J+1},\quad \phi_n\ket{2J+1}=0\quad(n>0),
\\
& L_0\ket{2J+1} = h_J \ket{2J+1}, \quad h_J = \frac{J(J+1)}{\kappa}.
\end{align*}
We also use abbreviated notations
$2J_{m,m'} + 1 = m - m'\kappa,\:
\F_{m,m'} = \F_{J_{m,m'}},\:
\Lambda_{m,m'} = \Lambda_{J_{m,m'}}, \:
L_{m,m'} = L_{\Lambda_{m,m'}}.$

Wakimoto representation is given by the following theorem  \cite{W2} (See also \cite{FeFr}).
\begin{thm}[Wakimoto representation]
Let $J,\kappa$ be complex numbers with $\kappa\neq 2,0$. Then the following $(\pi_\kappa, \F_{J})$ gives a level $\kappa-2$ representation of $\sltwohat$ on the Fock space $\F_{J}$:
\begin{align*}
& \pi_\kappa\bigl( e(\zeta) \bigr) = \sum_{n\in\Z} \pi_\kappa( e[n] ) \,\zeta^{-n-1} = \beta(\zeta),
\\
& \pi_\kappa\bigl( h(\zeta) \bigr) = \sum_{n\in\Z} \pi_\kappa( h[n] ) \,\zeta^{-n-1}
                  =-2:\gamma(\zeta)\beta(\zeta): + \kappa\partial\phi(\zeta),
\\
& \pi_\kappa( f(\zeta) \bigr) = \sum_{n\in\Z} \pi_\kappa( f[n] ) \,\zeta^{-n-1}
     =-:\gamma(\zeta)^2\beta(\zeta): + :\kappa\gamma(\zeta)\partial\phi(\zeta):
	   + (\kappa-2 )\partial\gamma(\zeta),
\\
& \pi_\kappa( d ) = L_0	
\end{align*}
where $:\quad :$ denotes the normal ordering.
This representation is isomorphic to the Verma module with highest weight
$\Lambda_J = (\kappa-2-2J)\Lambda_0 + 2J\Lambda_1- h_{J}\delta.$
\end{thm}

We define the screening current $S(\zeta)$ by
\begin{align*}
& S(\zeta)=\beta(\zeta):e^{-\phi(\zeta)}:
\end{align*}
The screening charge $Q_n$ can be defined as a map
$Q_n:\,\F_{m,m'}\mapsto \F_{m-2n,m'}$
 when $n\equiv m\mbox{ mod }p$ by
\begin{align*}
Q_n = \int_{\calC} \prod_{j=1}^{n} d\xi_j S(\xi_1) \cdots S(\xi_n).
\end{align*}
 and commutes with $\sltwohat$:\,
$[Q_n, x] = 0 \; (x \in\sltwohat).$
The contour $\calC$ is given in Figure \ref{fig:Felder}.
\begin{figure}[hbtp]
 \begin{center}
  \includegraphics[height=4cm]{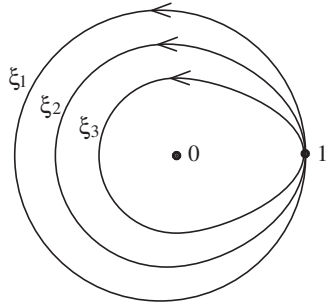}
    \caption{The contour $\calC$ for $k=4$.}
 \lb{fig:Felder}
 \end{center}
\end{figure}

To realize a irreducible representation with a dominant integral highest weight, we use the following theorem due to \cite{BeFe}.
\begin{thm}[BRST/BGG-resolution]
If the level is rational
\begin{align*}
&\kappa=p_1/p_2
\end{align*}
where $p_1,p_2$ are coprimes, then for $1\leq n\leq p_1-1,\, 0\leq n'\leq p_2-1$ the following sequence is exact
\begin{align*}
\begin{CD}
\cdots @>>> \overset{-1}{\F_{2p_1-n,n'}} @>{Q_{p_1-n}}>> \overset{0}{\F_{n,n'}} @>{Q_{n}}>> \overset{1}{\F_{-n,n'}}
@>>> \cdots
\end{CD}.
\end{align*}
And its zero-th cohomology group is isomorphic to the irreducible highest weight representation
\be
\mbox{{\rm H}}^i=
\begin{cases}
& L_{n,n'},\quad i=0,\\
& 0,\quad i\neq 0.
\end{cases}
\en
\end{thm}

Let us consider bosonic counter part of the intertwiner $\Phi$ of \eqref{eqn:Phi}.
This operator is called vertex operator and denoted by $\Phibar$. It satisfies
\begin{align*}
& \Phibar(\zeta): \F_{m,m'} \to \F_{m,m'} \otimes U(\zeta),
\\
& (x\otimes 1 + 1\otimes x)\, \Phibar(\zeta) = \Phibar(\zeta)\, x \quad (x\in\sltwohat'),
\end{align*}
and given by
\begin{align*}
& \Phibar(\zeta) = \sum_{n=0}^{2(k-1)} \Phibar_n(\zeta) \otimes u_n,
\\
& \Phibar_n(\zeta) = \frac{ (-)^{n} }{n!} \int_{\calC} \prod_{j=1}^{k-1} d\xi_j
                                        V_{n}(\zeta) S(\xi_{1})S(\xi_{2})\cdots S(\xi_{k-1}) ,
\\
& V_n(\zeta) = \gamma(\zeta)^{ 2(k-1) - n } :e^{(k-1)\phi(\zeta)}:.										
\end{align*}
In the above we took basis of $U(\zeta)$ as
\begin{align*}
& U(\zeta) = \bigoplus_{ n=0 }^{ 2(k-1) } \C u_n \otimes \C[\zeta,\zeta^{-1}],
\\
& u_n = f^n u_{0},\quad h u_0 = 2(k-1)\Lambdabar_1 u_0.
\end{align*}

\subsection{Integral Formula}
From now on we consider representations which satisfy the conditions
\begin{gather*}
2J+1=m \in\Z, \;m'=0,\;p_2=1,
\\
\kappa=p_1\in\Z_{\geq 0}
\end{gather*}
and abbreviate notations as $\F_m=\F_{m,m'=0},\: L_m= L_{m,m'=0}.$

\begin{prop}
\label{prop:intgrl}
We can naturally consider $\Phibar(\zeta)$ as an intertwiner for a irreducible highest weight representation:
\begin{align*}
\Phibar(\zeta)\in\Hom (L_{m,m'}, L_{m,m'} \otimes U(\zeta)).
\end{align*}
Then the trace function $\Jbar_{-m} (z|\tau)$ defined as
\begin{align}
\Jbar_{-m} (z|\tau) &= \tr_{L_m}\Bigl( \Phibar_{k-1}(1)\: p^{L_0-c_V/24}e^{2\pi iz h[0]} \Bigr)
\end{align}
is given by
\begin{align}
& \Jbar_{-m} (z) = I_{-m}(z|\tau) - e^{-2\pi i\frac{m(k-1)}{\kappa}} I_{m}(z|\tau),
\\
& I_{-m}(z|\tau) = \frac{1}{i\vartheta_1(2z)} \int_0^1\prod_{j=1}^{k-1}dv_j\;
	            \vartheta_{-m,\kappa} \Bigl( -2z+\frac 2\kappa V \Bigm|\tau \Bigr) F(v_1,\cdots,v_{k-1};z | \tau),
\label{eqn:I-integral}		
\displaybreak[2]		
\\[2mm]
& F(v_1,\cdots,v_{k-1} ;z| \tau) = (-)^{k-1}
\prod_{j=1}^{k-1}E(-v_j)^{-\frac{2(k-1)}{\kappa}}
                 \prod_{1\leq j<j'\leq k-1}E(v_j-v_{j'})^{\frac{2}{\kappa}}
				 \times \prod_{j=1}^{k-1}G(v_j;z|\tau),
\nonumber				
\\[2mm]
&
E(v)=2\pi i\frac{\vartheta_1(v)}{\vartheta_1'(0)},\quad
G(v;z|\tau)=
             \frac{\vartheta_1'(0)\vartheta_1(v+2z)}{\vartheta_1(2z)\vartheta_1(v)},
\nonumber
\\[2mm]
& \vartheta_{-m,\kappa} (x | \tau) = \sum_{ l\in\Z +\frac{m}{2\kappa} } p^{\kappa l^2} e^{2\pi
i(-x)\kappa l},
\quad V=\sum_{j=1}^{k-1} v_j.
\nonumber
\end{align}
The contour is the result of the change of variables: $\xi_j = e^{2\pi i v_j}$ and depicted in Figure \ref{fig:01}.
The definition of $\vartheta_1(v)=\vartheta_1(v|\tau)$ is given in Appendix \ref{appendix}. Differentiation of $\vartheta_1'(0)$ is taken w.r.t $v$.
\end{prop}
\begin{figure}[hbtp]
 \begin{center}
  \includegraphics[height=3.5cm]{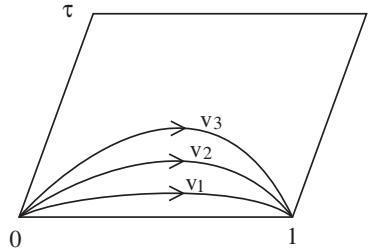}
  \caption{The contour for $k=4$.}
 \lb{fig:01}
 \end{center}
\end{figure}

\begin{proof}
\noi {\it Step 1}:\ We have
\begin{align*}
\tr_{\F_m} & \Bigl( V_{k-1}(\zeta) S(\xi_{1})S(\xi_{2})\cdots S(\xi_{k-1})
                                                                            p^{L_0}e^{2\pi iz
h[0]} \Bigr)
\\
= & \frac{p^{c_V/24}}{i\vartheta_1(2z)}
   e^{i\pi m(2z-\frac 2\kappa V)+i\pi\tau\frac{m^2}{2\kappa}}
   e^{2\pi iy\frac{(k-1)(m-k)}{\kappa}}
\\
 & \times\prod_{j=1}^{k-1}E(y-v_j)^{-\frac{2(k-1)}{\kappa}}
                 \prod_{1\leq j<j'\leq k-1}E(v_j-v_{j'})^{\frac{2}{\kappa}}
				 \times (k-1)!\prod_{j=1}^{k-1} \frac{1}{2\pi i\xi_j} G(v_j -y;z|\tau)
\end{align*}
where $\xi_j = e^{2\pi i v_j},\; \zeta = e^{2\pi i y}.$

\noi {\it Step 2}:\ Set
\begin{align*}
& \Xi(\zeta;z|\tau) = \Phibar_{k-1}(\zeta)\: p^{L_0-c/24}e^{2\pi iz h[0]},
\\
&
\Xi_{2l}(\zeta;z|\tau)=\Xi(\zeta;z|\tau),\quad \Xi_{2l-1}(\zeta;z|\tau)=e^{-2\pi i\frac{m(k-1)}{\kappa}} \Xi(\zeta;z|\tau).
\end{align*}
Then the following diagram is commutative
\[
\begin{CD}
\cdots @>{Q^{m}}>> \overset{-1}{\F_{2\kappa-m}}
  @>{Q^{\kappa-m}}>> \overset{0}{\F_{m}}
  @>{Q^{m}}>> \overset{1}{\F_{-m}}
  @>{Q^{\kappa-m}}>> \cdots
\\
@.                               @VV{\Xi_{-1}}V
                            @VV\Xi_{0}V
                        @VV\Xi_{1}V
                            @.
\\
\cdots @>{Q^{m}}>> \F_{2\kappa-m}
  @>{Q^{\kappa-m}}>> \F_{m}
  @>{Q^{m}}>> \F_{-m}
  @>{Q^{\kappa-m}}>> \cdots\;.			
\end{CD}
\]
Hence
\begin{align*}
\tr_{L_m}\Bigl(\Xi(\zeta;z|\tau)\Bigr) &= \sum_{l\in\Z} \tr_{\F_{m-2l\kappa}} \Bigl(\Xi_{2l}(\zeta;z|\tau)\Bigr)
                                                                    - \sum_{l\in\Z} \tr_{\F_{-m-2l\kappa}}\Bigl(\Xi_{2l-1}(\zeta;z|\tau)\Bigr).
\end{align*}
\end{proof}
Let
$\kappa = K+2k, \, K\in\Z_{\geq 0}, \,\lambda = (l-k) \bar{\Lambda}_1$.
\begin{thm}[Integral formula for $\Jhat_{\lambdahat}$]
\label{thm:intgrl}
\begin{align}
J_{\mu,K}(zh[0],\tau) &= \frac{\Jbar_{-m} (z|\tau)}{ \bra{m}\Phibar_{k-1}(1)\ket{m} } \frac{1}{\deltahat^{\prime k-1}}
\label{eqn:intgrl-1}
\end{align}
where $\Jbar_{-m} (z|\tau)$ is given by Prop.\ref{prop:intgrl} and
\begin{align}
\bra{m}\Phibar_{k-1}(1)\ket{m} &= \int_{\calC} \prod_{j=1}^{k-1} d\xi_j
               \prod_j^{k-1} \xi_j^{ \frac{-m+1}{\kappa} } (1-\xi_j)^{-\frac{2(k-1)}{\kappa} -1}
               \prod_{1\leq j<j'\leq k-1} (\xi_j - \xi_{j'})^{\frac{2}{\kappa}}.
\label{eqn:intgrl-2}
\end{align}
\end{thm}

\section{Modular transformation}
\label{sec:modular}
We study modular transformation property of $\Jbar_{-m} (z|\tau)$ using the technique proposed in \cite{FeSi}.

To simplify calculations we introduce $j_{-m}^{[s]}$ and $i_{-m}^{[s]}$
\begin{align*}
& j_{-m}^{[s]}(z|\tau) = i_{-m}^{[s]}(z|\tau) - i_{m}^{[s]}(z|\tau),
\\
& i_{-m+s}^{[s]}(z|\tau) = e^{-\pi is(-s+2m+\kappa-1)/2\kappa} e^{\pi im(k-1)/\kappa} \:I_{-m}^{[s]}(z|\tau)
\end{align*}
where
\begin{align*}
I_{-l}^{[s]}(z|\tau) = \frac{1}{i\vartheta_1(2z)}
                \int_0^\tau \prod_{j=k-s}^{k-1}dv_j \: \int_0^1\prod_{j=1}^{k-s-1}dv_j\;
	            \vartheta_{-m,\kappa} \Bigl( -2z+\frac 2\kappa V \Bigm| \tau \Bigr) F(v_1,\cdots,v_{k-1};z | \tau).
\end{align*}
Note that $j_{-m} (z|\tau) = e^{\pi im(k-1)/\kappa}\, \Jbar_{-m} (z|\tau).$

\subsection{$S$ transformation}
\begin{lem}
\label{lem:s-red}
\begin{align*}
2\sin\Bigl(\frac \pi \kappa l\Bigr)\: i_{-l}^{[s]}(z|\tau) &= i_{-l+1}^{[s-1]}(z|\tau) - i_{-l-1}^{[s-1]}(z|\tau),
\\[2mm]
2\sin\Bigl(\frac \pi \kappa l\Bigr)\: j_{-l}^{[s]}(z|\tau) &= j_{-l+1}^{[s-1]}(z|\tau) - j_{-l-1}^{[s-1]}(z|\tau).
\end{align*}
\begin{proof}
The lemma follows from
\begin{gather}
(1-e^{ \frac{2\pi i}{\kappa}(s-l) })\, I_{-l}^{[s]}(z|\tau) =
     I_{-l}^{[s-1]}(z|\tau) - e^{ \frac{2\pi i}{\kappa}(s-k) } I_{-l+2}^{[s-1]}(z|\tau).
\label{eqn:I-1}
\end{gather}
Consider contours for $v_{k-s}$ given in Figure \ref{fig:gamma}, \ref{fig:gamma2} and \ref{fig:gamma3}.
\begin{figure}[hbtp]
 \begin{minipage}{.50\linewidth}
  \begin{center}
   \includegraphics[width=.70\linewidth]{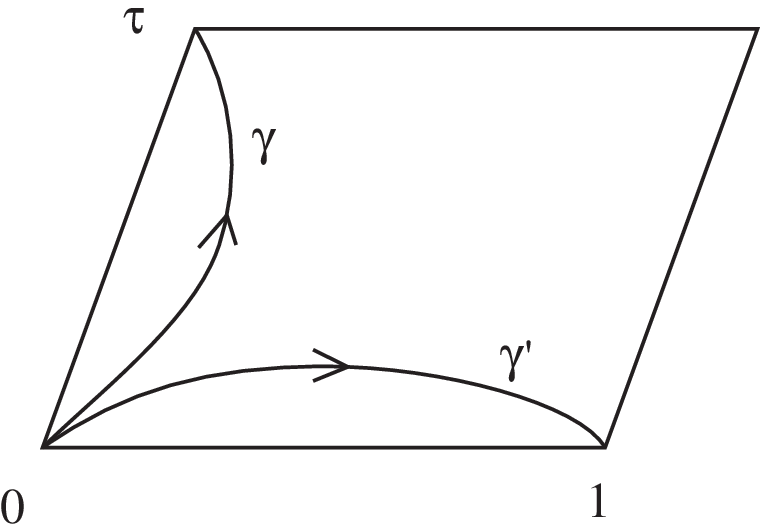}
   \caption{The contours $\gamma$ and $\gamma'$.}
   \lb{fig:gamma}
  \end{center}
 \end{minipage}
 \begin{minipage}{.50\linewidth}
  \begin{center}
   \includegraphics[width=.70\linewidth]{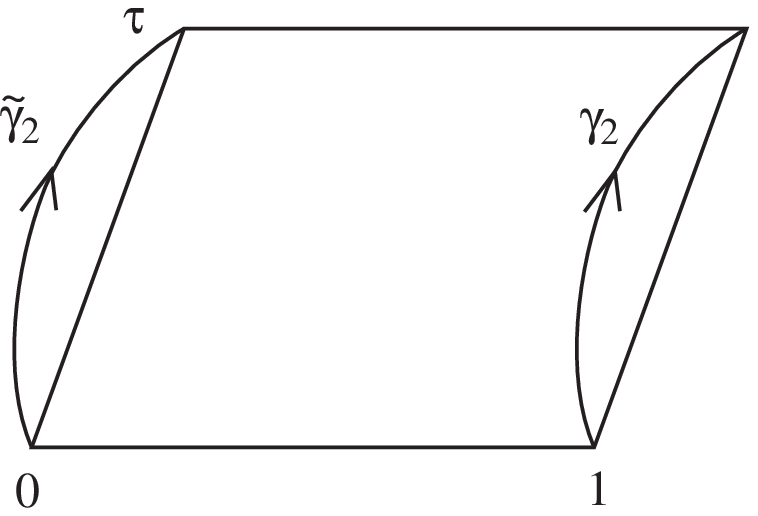}
   \caption{The contours $\gamma_2$ and $\widetilde{\gamma_2}$.}
   \lb{fig:gamma2}
  \end{center}
 \end{minipage}
\end{figure}

\begin{figure}
 \begin{center}
  \includegraphics[height=3.5cm]{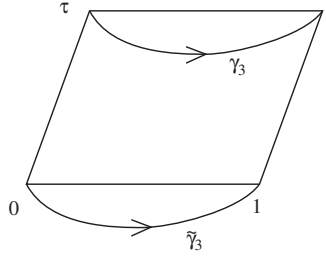}
  \caption{The contours $\gamma_3$ and $\widetilde{\gamma_3}$.}
 \lb{fig:gamma3}
 \end{center}
\end{figure}

Contours for $I_{-l}^{[s]}(z|\tau)$ is as follows. The variable $v_{k-s}$ is on $\gamma$. The contours of $v_j\: (k-s+1\leq j\leq k-1)$ are between $[0,\tau]$ and $\gamma$ from right to left. The contours of $v_j\: (1\leq j\leq k-s-1)$ are between $[0,1]$ and $\gamma'$ from bottom to top.

We deform the contour of $v_{k-s}$ in $I_{-l}^{[s]}(z|\tau)$ as
\begin{align}
& \int\! \prod_{j\neq k-s}\! dv_j
\int_{\gamma} dv_{k-s}\,
 	\vartheta_{-l,\kappa} \Bigl( -2z+\frac 2\kappa V |\tau \Bigr)
	F(\cdots,v_{k-s},\cdots;z|\tau)
\nonumber
\\[2mm]
& = \int\! \prod_{j\neq k-s}\! dv_j
    \left( \int_{\gamma'} + \int_{\gamma_2} - \int_{\gamma_3} \right) dv_{k-s}\,
	\vartheta_{-l,\kappa} \Bigl( -2z+\frac 2\kappa V \Bigm| \tau \Bigr)
\nonumber	
\\[2mm]	
& \hspace{7cm}	 
	\times  F(\cdots,v_{k-s},\cdots;z|\tau)
\label{eqn:I-2}
\end{align}
Then by carefully studying crossing of contours, we can show
\begin{align*}
& \int\! \prod_{j\neq k-s}\! dv_j
    \int_{\gamma_2} dv_{k-s}\,
	\vartheta_{-l,\kappa} \Bigl( -2z+\frac 2\kappa V \Bigm| \tau \Bigr)
	F(\cdots,v_{k-s},\cdots;z|\tau)
\\[2mm]
&= \int\! \prod_{j\neq k-s}\! dv_j
       \int_{\widetilde{\gamma_2}} dv_{k-s}\,
	   \vartheta_{-l,\kappa} \Bigl( -2z+\frac 2\kappa (V+1) \Bigm| \tau \Bigr)
       F(\cdots,v_{k-s}+1,\cdots;z|\tau)
\\[2mm]	
&= e^{ -\frac{2\pi i}{\kappa}(s-1) }
	\int\! \prod_{j\neq k-s}\! dv_j
	\int_{\gamma} dv_{k-s}\,
	e^{-2\pi i \frac l\kappa} \vartheta_{-l,\kappa} \Bigl( -2z+\frac 2\kappa V \Bigm| \tau \Bigr)
\\[2mm]	
& \hspace{7cm}	 
	\times  e^{ \frac{2\pi i}{\kappa}(2s-1) }  F(\cdots,v_{k-s},\cdots;z|\tau),
\end{align*}
and
\begin{align*}
& \int\! \prod_{j\neq k-s}\! dv_j
	\int_{\gamma_3} dv_{k-s}\,
    \vartheta_{-l,\kappa} \Bigl( -2z+\frac 2\kappa V \Bigm| \tau \Bigr)
	F(\cdots,v_{k-s},\cdots;z|\tau)
\\[2mm]
&= \int\! \prod_{j\neq k-s}\! dv_j
		\int_{\widetilde{\gamma_3}} dv_{k-s}\,
      \vartheta_{-l,\kappa} \Bigl( -2z+\frac 2\kappa (V+\tau) \Bigm| \tau \Bigr)
      F(\cdots,v_{k-s}+\tau,\cdots;z|\tau)
\\[2mm]
&= e^{ -\frac{2\pi i}{\kappa}(k-s) }
		\int\! \prod_{j\neq k-s}\! dv_j
		\int_{\gamma'} dv_{k-s}\,
    e^{ -2\pi i\frac{\tau}{\kappa} } e^{ -2\pi i (-2z+\frac{2}{\kappa}V) }
	    \vartheta_{-l+2,\kappa} \Bigl( -2z+\frac 2\kappa V \Bigm| \tau \Bigr)
\\
	&	\hspace{6cm} \times
    e^{ \frac{2\pi i}{\kappa}(2V+\tau-1) } e^{-4\pi iz} F(\cdots,v_{k-s},\cdots;z|\tau).
\end{align*}
Throughout this calculation contours of $v_j\:(j\neq k-s)$ remain the same as in $I_{-l}^{[s]}(z|\tau)$.
Plugging in these equations to \eqref{eqn:I-2}, we obtain \eqref{eqn:I-1}.
\end{proof}
\end{lem}

\begin{lem}
\label{lem:i-mod}
\begin{align*}
i_{-m} \Bigl( \frac z\tau \Big| -\frac 1\tau \Bigr) &=
   i\, e^{ \frac{\pi i}{2} \frac{(k-1)(\kappa-k)}{\kappa} }
  \tau^{\frac{(k-1) k}{\kappa}}
   e^{2\pi i\frac{\kappa-2}{\tau} z^2} \frac{1}{\sqrt{2\kappa}}
   \sum_{l=0}^{2\kappa-1} e^{-\frac{\pi i}{\kappa} ml} \: i_{-l}^{[k-1]}(z|\tau),
\\
i_{m} \Bigl( \frac z\tau \Big| -\frac 1\tau \Bigr) &=
   i\, e^{ \frac{\pi i}{2} \frac{(k-1)(\kappa-k)}{\kappa} }
  \tau^{\frac{(k-1) k}{\kappa}}
   e^{2\pi i\frac{\kappa-2}{\tau} z^2} \frac{1}{\sqrt{2\kappa}}
   \sum_{l=0}^{2\kappa-1} e^{-\frac{\pi i}{\kappa} ml} \: i_{l}^{[k-1]}(z|\tau),
\end{align*}
\end{lem}

\begin{proof}
Transformation of theta functions under the action of $S$ is given by
\begin{align*}
& \vartheta_1 \Bigl( \frac v\tau \Big| - \frac 1\tau \Bigr)
      = \frac 1i \sqrt{\frac \tau i} e^{\frac{i\pi}{\tau} v^2} \vartheta_1 (v|\tau),
\\[2mm]
& \vartheta_1' \Bigl( 0\Big| - \frac 1\tau \Bigr)	
      = \frac\tau i \sqrt{\frac \tau i} \vartheta_1' (0|\tau),
\\[2mm]
& E	\Bigl( -\frac v\tau \Big| - \frac 1\tau \Bigr)
     =\tau^{-1} e^{\frac{i\pi}{\tau} v^2} E(-v|\tau),
\\[2mm]
& G \Bigl( \frac v\tau ; \frac z\tau\Big| - \frac 1\tau \Bigr)
      = \tau e^{\frac{\pi i}{\tau} 4vz} \: G(v:z|\tau),
\\[2mm]
& \vartheta_{-m,\kappa} \Bigl( \frac x\tau \Big| -\frac 1\tau \Bigr)
   = \sqrt{\frac \tau i}\, e^{i\pi \frac{\kappa}{2\tau} x^2}
\frac{1}{\sqrt{2\kappa}}
      \sum_{l=0}^{2\kappa-1} e^{-\frac{\pi i}{\kappa} ml} \: \vartheta_{-l,\kappa} (x|\tau) .
\end{align*}
Here the branch of $\sqrt{\tau/ i}$ is $\sqrt{\tau/ i}=1\mbox{ for }\tau=i.$
Hence from \eqref{eqn:I-integral}, with the change of a variable $v_j\mapsto v_j/\tau$, we obtain
\begin{align*}
I_{-m} \Bigl( \frac z\tau \Big| -\frac 1\tau \Bigr) =
   i\: \tau^{\frac{(k-1) k}{\kappa}} e^{ 2\pi i \frac{\kappa-2}{\tau} z^2 } \frac{1}{\sqrt{2\kappa}}
   \sum_{l=0}^{2\kappa-1} e^{-\frac{\pi i}{\kappa} ml} \: I_{-l}^{[k-1]}(z|\tau)
\end{align*}
where
\[
\tau^{\frac{(k-1) k}{\kappa}}=
       \underbrace{ \tau^{\frac{2(k-1)}{\kappa}}\cdots\tau^{\frac{2(k-1)}{\kappa}} }_{k-1}
       \underbrace{ \tau^{-\frac{2}{\kappa}}\cdots\tau^{-\frac{2}{\kappa}} }_{\binom{k-1}{2}}
\]
and $\tau^{\frac{2(k-1)}{\kappa}},\;\tau^{-\frac{2}{\kappa}}$ are of certain  branches.

For $i_{-m}(z|\tau)$, we use
$i_{-m+2\kappa}(z|\tau)=i_{-m}(z|\tau)$.
\end{proof}
From these lemmas, this proposition follows.
\begin{prop}
\label{prop:jS}
\begin{align}
j_{-m} \Bigl( \frac z\tau \Big| -\frac 1\tau \Bigr) &=
   i\, e^{ \frac{\pi i}{2} \frac{(k-1)(\kappa-k)}{\kappa} }
  \tau^{\frac{(k-1) k}{\kappa}}
   e^{ 2\pi i \frac{\kappa-2}{\tau} z^2 } \frac{1}{\sqrt{2\kappa}}
   \sum_{l=0}^{2\kappa-1} \Bigl(UA^{k-1}\Bigr)_{m,l} \: j_{-l}(z|\tau),
\label{eqn:jS}
\\
U_{m,l} &= e^{-\frac{\pi i}{\kappa} ml},
\nonumber
\\
A_{m,l} &= \delta_{m-1,l} \Bigl( 2\sin\Bigl(\frac \pi \kappa m \Bigr) \Bigr)^{-1}
                   - \delta_{m+1,l} \Bigl( 2\sin\Bigl(\frac \pi \kappa m \Bigr) \Bigr)^{-1},
\nonumber
\end{align}
where $\delta_{-1,2\kappa-1}=\delta_{2\kappa,0}=1$.
\end{prop}
To rewrite the sum in the R.H.S. of \eqref{eqn:jS} as a sum over $\lambda\in\Phat^+_K$, we need one more proposition.
\begin{prop}
\label{prop:fusion}
\begin{itemize}
\item[]
\item[(1)]
  Fusion rule:
  \begin{align*}
  &  j_{-m}^{[s]} (z|\tau) = 0 \quad\mbox{  for  }\quad 0\leq m\leq (k-1)-s,
  \\
  &  j_{-(\kappa-m)}^{[s]} (z|\tau) = 0 \quad\mbox{  for  }\quad 0\leq m\leq (k-1)-s.
  \end{align*}
\item[(2)]
$j_{m+\kappa}^{[s]} (z|\tau) = -j_{\kappa-m}^{[s]} (z|\tau).$
\end{itemize}

\begin{proof}
(1)
We prove an equivalent statements:
\begin{align}
& i_{-m}^{[s]} (z|\tau) = i_{m}^{[s]} (z|\tau) \quad\mbox{for}\quad 0\leq m\leq (k-1)-s,
\label{eqn:fusion-pf1}\\
& i_{-\kappa-m}^{[s]} (z|\tau) = i_{-\kappa+m}^{[s]} (z|\tau)
   \quad\mbox{for}\quad 0\leq m\leq (k-1)-s
\label{eqn:fusion-pf2}
\end{align}
by induction.
\\
(i) $l=0$ trivial.
\\
(ii) $l=1$ Set $l=0$ in Lemma \ref{lem:s-red}.
\\
(iii) Suppose the statement is true for $l\leq m$, then
\begin{align*}
& i_{-m+1}^{[s-1]}(z|\tau) - i_{-m-1}^{[s-1]} (z|\tau)
\\[2mm]
& = 2\sin\Bigl(\frac\pi \kappa m \Bigr)\: i_{-m}^{[s]}(z|\tau)
\\[2mm]
& = 2\sin\Bigl(\frac\pi \kappa m\Bigr)\: i_{m}^{[s]}(z|\tau)
\\[2mm]
&= -i_{m+1}^{[s-1]}(z|\tau) + i_{m-1}^{[s-1]} (z|\tau)
\end{align*}
the first and third equality is valid for $1\leq s\leq k-1$. Hence
\begin{align*}
& i_{-m+1}^{[s-1]}(z|\tau) - i_{-m-1}^{[s-1]} (z|\tau)
= -i_{m+1}^{[s-1]}(z|\tau) + i_{m-1}^{[s-1]} (z|\tau)
\quad\mbox{  for  }\quad 1\leq s\leq (k-1)-m.
\end{align*}
Combining this with
\[
i_{-m-1}^{[s-1]} (z|\tau) = i_{m+1}^{[s-1]} (z|\tau) \quad\mbox{  for  }\quad 1\leq s\leq
(k-1)-m,
\]
we get
\[
i_{-m}^{[s]} (z|\tau) = i_{m}^{[s]} (z|\tau) \quad 0\leq s\leq (k-1)-m.
\]
which is equivalent to \eqref{eqn:fusion-pf1}. A similar argument with $l=\kappa$ instead of $0$ in Lemma \ref{lem:s-red} gives \eqref{eqn:fusion-pf2}.
\\
(2)
\begin{align*}
j_{m+\kappa}^{[s]} (z|\tau)
=& \,i_{m+\kappa}^{[s]} (z|\tau) - i_{-m-\kappa}^{[s]} (z|\tau)
\\
=& \,i_{m-\kappa}^{[s]}(z|\tau) - i_{-m+\kappa}^{[s]} (z|\tau)
\\
=& -j_{\kappa-m}^{[s]} (z|\tau).
\end{align*}
\end{proof}
\end{prop}
From Proposition \ref{prop:jS} and \ref{prop:fusion}, we obtain the main statement of this subsection.
\begin{prop}
\label{prop:S-1}
\begin{align}
j_{-m} \Bigl( \frac z\tau \Big| -\frac 1\tau \Bigr) &=
   \, e^{ \frac{\pi i}{2} \frac{(k-1)(\kappa-k)}{\kappa} }
  \tau^{\frac{(k-1) k}{\kappa}}
   e^{ 2\pi i \frac{\kappa-2}{\tau} z^2 } \frac{2}{\sqrt{2\kappa}}
   \sum_{l=k}^{\kappa-k} S(K,k)_{m,l} \: j_{-l}(z|\tau),
\\[2mm]
S(K,k) &= U(K,k) A(K,1) A(K,2) \cdots A(K,k-1),
\\[2mm]
U(K,k) &=
\begin{pmatrix}
s(k)                 &    s(2k)    & \hdotsfor{2} & s(k(\kappa-1)) \\
s(k+1)            & \ddots    & \hdotsfor{2} & s((k+1)(\kappa-1)) \\
\vdots           &                  & \ddots   &     & \vdots    \\
s(\kappa-k)  &            \hdotsfor{3}          & s((\kappa-k)(\kappa-1))
\end{pmatrix},
\displaybreak[3]
\\[2mm]
A(K,l) &= \frac 12
\begin{pmatrix}
-s(l)^{-1}    &                        &                  &        &               \\
0                  &  -s(l+1)^{-1}  &                 &        &               \\
s(l+2)^{-1} &  0                    & \ddots     &        &                \\
0                  & s(l+3)^{-1}     &  \ddots   &       & -s(\kappa-l-2)^{-1} \\
                    &                         & \ddots    &        &  0           \\
		            &                         &                  &        & s(\kappa-l)^{-1}     \\
\end{pmatrix}
\end{align}
where $s(k) = \sin( \pi k/\kappa ).$
\end{prop}

\noindent {\it Example 1} \quad
$k=1$: $S(K,1)=U(K,1)$ recovers the known result for the character.

\noindent {\it Example 2} \quad
$p=2k$: Using
\begin{align*}
S(2k,k) = (-1)^{k-1} \sqrt{k},
\end{align*}
we can be show that behavior of $j_{-m}(z|\tau)$ and $\hat{\delta}^{'k-1}$ under the action of $\Gamma$ are the same.

\vspace{5mm}
It is known that elements of the matrix for $S$-transformation are related to special values of Macdonald's polynomial \cite{Ki}\cite{FSV}. Here we give explicit expression for this interesting connection.

We adopt the definition and notation for Macdonald's polynomial $P_{\lambda}(x; q, t)$ in \cite{M}. Special values we consider is denoted by $P^{(k)}_l(m)$ and defined as follows:
\begin{align}
P^{(k)}_l(m) = P_{l\Lambdabar_1}(x_1=q^{ -m}, x_2=q^{m}; q^2, q^{2k}),\quad q=e^{\pi i/\kappa}.
\end{align}

\begin{prop}
\label{prop:S-2}
\begin{align}
S(K,k) &= (-2)^{k-1}
\begin{pmatrix}
\prod_{j=1}^{k-1} s(k-j)  &                                                    &                 &                                                 \\
                                                &  \prod_{j=1}^{k-1} s(k+1-j)    &                 &                                                 \\
                                                &                                                    &  \ddots   &                                                 \\
                                                &  		                                     	     &                  & \prod_{j=1}^{k-1} s(\kappa -k-j)
\end{pmatrix}
\\[2mm]
& \times
\begin{pmatrix}
  P^{(k)}_0(k)                      &  P^{(k)}_0(k+1)  &  \cdots  &  P^{(k)}_0(\kappa-k)                      \\
  P^{(k)}_1(k)                      &  \ddots               &                 &                                                             \\
  \vdots                               &                              &  \ddots  &                                                             \\
  P^{(k)}_{\kappa-2k}(k)  &                               &                &  P^{(k)}_{\kappa-2k}(\kappa-k)
\end{pmatrix}
\begin{pmatrix}
s(k) &                &                &                           \\
        &  s(k+1)  &                 &                           \\
        &                 & \ddots  &                            \\
        &  			     &               & s(\kappa-k)
\end{pmatrix}.
\nonumber
\end{align}
\end{prop}
Before giving proof we prepare two lemmas. Set
\begin{align*}
S(K,k) &=
\begin{pmatrix}
  \sigma^{(k)}_{k,k}                 &   \sigma^{(k)}_{k,k+1}    &    \cdots      &   \sigma^{(k)}_{k,\kappa-k}                \\
  \sigma^{(k)}_{k+1,k}             &   \ddots                            &    & \\
  \vdots                                     &                                             &  \ddots
&                                                                   \\
  \sigma^{(k)}_{\kappa-k,k}  &   &   &
\sigma^{(k)}_{\kappa-k,\kappa-k}
\end{pmatrix}.
\end{align*}
\begin{lem}
\begin{align*}
  \sigma^{(k)}_{n,m} = \frac{\sigma^{(k-1)}_{n,m+1}}{2s(m+1)} - \frac{\sigma^{(k-1)}_{n,m-1}}{2s(m-1)}.
\end{align*}
\end{lem}
\begin{proof}
Use
\begin{align*}
S(K,k-1) A(K,k-1) =
\left(
\begin{array}{c}
  * \quad\cdots\cdots\quad *
\\
  \boxed{
  \begin{array}{c}
       \\ S(K, k+1) \\ {}
  \end{array}
  }
\\
  * \quad\cdots\cdots\quad *
\end{array}
\right).
\end{align*}
\end{proof}
\begin{lem}
\begin{align*}
  P^{(k)}_{n}(x+1) - P^{(k)}_{n}(x-1) = 4\,s(-n)s(x) P^{(k+1)}_{n-1}(x).
\end{align*}
\end{lem}
\begin{proof}
Theorem 4.1 of \cite{FSV}.
\end{proof}
\begin{proof}
We prove by induction
\begin{align*}
   \sigma^{(k)}_{n,m} = 2^{k-1} s(m) \prod_{j=1}^{k-1} s(-n+j) \times  P^{(k)}_{n-k}(m).
\end{align*}
\noindent {\it k=1}:
\begin{align*}
\mbox{L.H.S} = s(nm) =s(m) \frac{s(nm)}{s(m)} = s(m) P^{(1)}_{n-1}(m) = \mbox{R.H.S}.
\end{align*}
\noindent {\it Step of induction}:
\begin{align*}
\sigma^{(k)}_{n,m} & = \frac{\sigma^{(k-1)}_{n,m+1}}{2s(m+1)} - \frac{\sigma^{(k-1)}_{n,m-1}}{2s(m-1)}
\\
& = 2^{k-3} \prod_{j=1}^{k-2} s(-n+j) \times  \bigl( P^{(k-1)}_{n-(k-1)}(m+1) - P^{(k-1)}_{n-(k-1)}(m-1) \bigr)
\\
& = 2^{k-1} s(m) \prod_{j=1}^{k-1} s(-n+j) \times  P^{(k)}_{n-k}(m).
\end{align*}
\end{proof}
\subsection{$T$ transformation}
\begin{prop}
\label{prop:jT}
\begin{align*}
j_{-m}(z|\tau+1) =& e^{ \pi i \bigl(\frac{m^2}{2\kappa} - \frac{1}{4}\bigr) } j_{-m}(z|\tau)
\end{align*}
\end{prop}
\begin{proof}
Use
\begin{align*}
\vartheta_1 (v|\tau+1) &= e^{\pi i/4} \vartheta_1 (v|\tau),
\\
\vartheta_{m,\kappa} (v|\tau+1) &= e^{\pi i \frac{m^2}{2\kappa}} \vartheta_{m,\kappa} (v|\tau).
\end{align*}
\end{proof}

\subsection{Affine Jack's polynomials}
We rewrite the integral formula \eqref{eqn:intgrl-1} \eqref{eqn:intgrl-2} using $j_{-m} (z|\tau)$:
\begin{gather}
J_{\lambda, K} (h,u,\tau) = e^{2\pi i(\kappa-2k)u} J_{\lambda, K} (h,\tau),
\nonumber
\\
J_{\lambda, K} (h,\tau) = g_{l,K,k}\, \frac{ j_{-l}(z|\tau) }{ \hat{\delta}^{'k-1}},
\label{eqn:Jj}
\\
g_{l,K,k} = \frac{e^{-\pi il(k-1)/\kappa}}{ \bra{l}\Phibar_{k-1}(1)\ket{l} }
\label{eqn:g-def}
\end{gather}
where $\lambda = (l-k)\Lambdabar_1$.
From \eqref{eqn:intgrl-2}, if we take branches of the integrand appropriately, the factor $g_{m,K,k}$ is given by the following formula:
\begin{align}
g_{m,K,k}^{-1} =
  e^{\pi i (k-1) (\frac 1\kappa +1)} & (2i)^{k-1}
  \prod_{n=0}^{k-2} \sin\left( \pi\biggl( \frac{-m+1+n}{\kappa} + 1 \biggr)\right)
  \label{eqn:g}
\\
  & \times
  \prod_{n'=1}^{k-2} \sum_{j=0}^{n'} e^{\pi i \frac{n'-2j}{\kappa}}
  \:\times B_{k-1}\left( \frac{-m+1}{\kappa} + 1,-\frac{2(k-1)}{\kappa}, \frac{1}{\kappa}\right)
\nonumber
\end{align}
where
\begin{align}
B_n(\alpha,\beta,\gamma) & =
  \int_{\Delta_n} \prod_{j=1}^n t_j^{\alpha-1} (1-t_j)^{\beta-1} \prod_{0\leq j< k \leq n} (t_j-t_k)^{2\gamma}
\\
  & = \frac{1}{n!} \prod_{j=0}^{n-1}
  \frac{
    \Gamma(1+(1+j)\gamma) \Gamma(\alpha+j\gamma) \Gamma(\beta+j\gamma)
	}{
	\Gamma(1+\gamma) \Gamma(\alpha+\beta+(n+j-1)\gamma)
	}
\end{align}
and $\Delta_n = \{ t\in\R^n \,|\, 0\leq t_n < \cdots < t_1 \leq 1\}$.

We consider \eqref{eqn:intgrl-2} as analytic continuation of parameters $m,k,\kappa$ from the region where each curve of the contour $\calC$ can be deformed to $[0,1] + [1,0]$ and related to $\Delta_n$.
Completely rigorous treatment of this subject using the cohomology with coefficients in local systems is beyond the scope of this paper. Let us consider $k=2$ case as an example.
In this case we consider
\begin{align*}
\int d\xi\, \xi^{\alpha-1} (1-\xi)^{\beta-1}
\end{align*}
for arbitrary $\alpha,\beta$ as analytic continuation from $\alpha,\beta>0$.
If we take the Pochammer loop as the contour, it can be shrunk around $\xi=1$ since $\beta>0$ and gives the double of $\calC$.

Now we can state our main result about modular transformation of the affine Jack's polynomials.
\begin{thm}
\begin{itemize}
\item[]
\item[(1)]
$J_{\lambda, K} (z,u,\tau+1) = e^{ 2\pi i \frac{-k\kappa+2l^2}{8\kappa} } J_{\lambda, K} (z,u,\tau),$
\item[(2)]
Define $S^J(K,k)_{m,l}$ by
\begin{align*}
J_{\mu, K}  & \Bigl( \frac z\tau ,u-\frac{z^2}{\tau}, -\frac 1\tau \Bigr)
= \sum_{\lambda\in P^+_K} S^J(K,k)_{m,l} \: J_{\lambda, K} (z ,u, \tau)
\end{align*}
where $\lambda = (l-k)\Lambdabar_1,\,\mu = (m-k)\Lambdabar_1$. Then
\begin{align*}
S^J(K,k) &= (-1)^{k-1}
   \, e^{ \frac{\pi i}{2} (k-1)\frac{K}{2\kappa} }
   \tau^{ -(k-1)\frac{K}{2\kappa} }
   \frac{2}{\sqrt{2\kappa}}
   \left(G S(K,k) G^{-1}\right),
\\
G_{m,l} &= g_{m,K,k}\; \delta_{m,l}.
\end{align*}
\end{itemize}
\end{thm}
The point is that we can give an explicit expression for matrix elements:
\begin{align*}
\left(G S(K,k) G^{-1}\right)_{m,l} = \frac{g_{m,K,k}}{g_{l,K,k}} \, S(K,k)_{m,l}
\end{align*}
where $S(K,k)_{m,l}$ is given in Proposition \ref{prop:S-1} or \ref{prop:S-2} and
\begin{gather}
\frac{g_{m,K,k}}{g_{m+n,K,k}} =
  \frac{
  \displaystyle{
  \prod_{j=m+1-k}^{m+n-k} \Gamma\left( \frac{j}{\kappa} \right) \Gamma\left( \frac{K+k-j}{\kappa} \right)
  }}{
 \displaystyle{
  \prod_{j'=0}^{n-1} \Gamma\left( \frac{m+j'}{\kappa} \right) \Gamma\left( \frac{K+k-m-j'}{\kappa} \right)
  }}
\end{gather}
for $n\in\N$.
\begin{proof}
(1) Proposition \ref{prop:jT}.

(2) Based on Proposition \ref{prop:S-1} and \ref{prop:S-2}.
Note $\sum_{\lambda\in P^+_K} = \sum_{l=k}^{\kappa-k}$.

About the choice of a branch in the calculation, we note the following: Since $\hat{\delta}' = i\vartheta(2z|\tau)$
\begin{align*}
\left( \frac{j_{-m}}{\hat{\delta}^{' k-1}} \right) & \Bigl( \frac z\tau \Big| -\frac 1\tau \Bigr)
\\
&= \left( \frac{1}{i} \sqrt{\frac{\tau}{i}} \right)^{-(k-1)}
   \, e^{ \frac{\pi i}{2} \frac{(k-1)(\kappa-k)}{\kappa} }
   \,\tau^{\frac{(k-1) k}{\kappa}}
   e^{ 2\pi i \frac{\kappa-2k}{\tau} z^2 } \frac{2}{\sqrt{2\kappa}}
   \sum_{l=k}^{\kappa-k} S(K,k)_{m,l} \: \frac{j_{-l}}{\hat{\delta}^{' k-1}}(z|\tau).
\end{align*}
From the fact at $K=0$:
\begin{align*}
\frac{j_{-l}}{\hat{\delta}^{' k-1}}(z|\tau) =1,
\end{align*}
we have
\begin{align*}
\left( \frac{1}{i} \sqrt{\frac{\tau}{i}} \right)^{k-1}
= (-1)^{k-1} \tau^{\frac{k-1}{2}}  e^{ \frac{\pi i}{2} \frac{k-1}{2} }.
\end{align*}
Taking the factor $e^{2\pi i(\kappa-2k)u}$ into account, we obtain the formula.
\end{proof}

\subsection{Relation to \cite{FSV}}
Let us make some comments on the relation between the results of \cite{FSV} and ours.
We represent the objects of \cite{FSV} with the quotation marks \lq\ \rq.

In \cite{FSV} they studied holomorphic solutions of the KZB-heat equation with certain properties which are called conformal blocks.
A set of solutions \lq$u^{[p]}_n(\lambda,\tau)$\rq\  are given in terms of elliptic Selberg type
integrals. These were obtained in \cite{FV} from geometrical point of view.
It is shown that the solutions \lq$u^{[p]}_n(\lambda,\tau)$\rq\ form a basis of the space of conformal blocks which are invariant under the action of $SL(2,\Z)$.
And the projective action of the modular group with respect to this basis is computed and the relation to special values of Macdonald's polynomials is studied.

These results are completely parallel to ours obtained from representation theoretical viewpoint. Details of the correspondence are as follows.
For parameters and arguments, we have $\lq\kappa\mbox{\rq}=\kappa,\, \lq\tau\mbox{\rq}=\tau,\,\lq p\mbox{\rq}=k-1,\,\lq\lambda\mbox{\rq} =-2z$.
The KZB-heat equation is the eigenvalue problem for the operator $\Lhat_k$ of \eqref{eqn:L} instead of $\Mhat_k$.
For fixed $\kappa,\lq p\mbox{\rq}$, the space of conformal blocks and its basis
$\lq\{ u^{[p]}_n(\lambda,\tau) | p+1 \leq n \leq \kappa-p-1\}\mbox{\rq}$
correspond to the space of $W$-symmetric theta functions $A^{\What}_K$ and
its basis $\{ J_{\lambda,K} | \lambda\in P^+_K \}$.
The correspondence is in some sense the most apparent at the level of integrals:
$\lq J^{[k]}_{\kappa,n}(\lambda,\tau)\mbox{\rq} $ corresponds to
$I^{[s]}_{-l}(z|\tau)$ with $\lq k\mbox{\rq}=k-s-1,\,\lq n\mbox{\rq}=l$.

As for the action of $SL(2,\Z)$, the projective action on the space of conformal blocks is equivalent to the action described is Theorem 4.10.
The transition matrices are concretely given in Prop. 6.2 of \cite{FSV}.
(But in the paper \cite{Ki}, on which the calculation of Prop.6.2 is based,
the normalization of affine Jack's polynomials are not rigorously considered due to the limit of the technique.
Hence the formulas in Prop. 6.2 needs some corrections.)

\section{Affine Jack's polynomials in terms of modular and elliptic functions}
\label{sec:example}
It seems very difficult to perform integrals in the integral formula in Theorem \ref{thm:intgrl} since $\Jbar_{-m} (z|\tau)$ is given by multi-integral of multivalued function. But when the level $K$ is small, we can use modular transformation property to write down $J_{\lambda,K}$ in terms of modular and elliptic functions.

\subsection{Affine Jack's polynomials in terms of modular and elliptic functions}
The main result of this section is the following theorem which gives $J_{\lambda,K}$ in terms of characters of $\sltwohat$. Together with the integral formula given in Section \ref{sec:intgrl}, this gives some formulas for elliptic Selberg type integral. The proof is given in the subsequent subsections.
\begin{thm}
\label{thm:K=1,2}
Let\enskip $\chi_{\lambdahat}(z,u,\tau) = \tr_{L_{\lambdahat}}\Bigl( e^{2\pi i zh[0]} p^{L_0-\frac{c_V}{24}} \Bigr)
e^{2\pi i Ku}$ be the normalized character of $L_{\lambdahat}$.
\begin{itemize}
\item[(1)]
Level one:
\begin{align}
& J_{0, 1}(z,u,\tau) = \eta(\tau)^{-\frac{k-1}{\kappa}} \chi_{\Lambda_0}(z,u,\tau),
\\
& J_{\bar{\Lambda}_1, 1}(z,u,\tau)
= \eta(\tau)^{-\frac{k-1}{\kappa}} \chi_{\Lambda_1}(z,u,\tau).
\end{align}
Here $\eta(\tau) = p^{\frac{1}{24}} \prod_{j=0}^{\infty} (1-p^j)$ is the Dedekind's $\eta$-function.
\item[(2)]
Level two:
\begin{align}
& J_{0, 2}(z,u,\tau) =
\frac{
          \chi_{2\Lambda_0}(z,u,\tau) + \chi_{2\Lambda_1}(z,u,\tau)
}{2h_1(\tau)}
+
\frac{
          \chi_{2\Lambda_0}(z,u,\tau) - \chi_{2\Lambda_1}(z,u,\tau)
}{2h_2(\tau)},
\\[2mm]
& J_{\bar{\Lambda}_1, 2}(z,u,\tau)
= \left( \frac{1}{\sqrt{2}} \right)^{\frac{k-1}{k+1}}
      \frac{\chi_{\Lambda_0+\Lambda_1}(z,u,\tau)}{h_3(\tau)},
\\[2mm]
& J_{2\bar{\Lambda}_1, 2}(z,u,\tau) =
\frac{
          \chi_{2\Lambda_0}(z,u,\tau) + \chi_{2\Lambda_1}(z,u,\tau)
}{2h_1(\tau)}
-
\frac{
          \chi_{2\Lambda_0}(z,u,\tau) - \chi_{2\Lambda_1}(z,u,\tau)
}{2h_2(\tau)},
\intertext{where}
&
h_1(\tau) = \dstyle{ \left( \frac{\eta(\tau/2)\eta(2\tau)}{\eta(\tau)}\right)^{\frac{k-1}{k+1}}, }\quad
h_2(\tau) = \dstyle{ \left( \frac{\eta(\tau)^2}{\eta(\tau/2)} \right)^{\frac{k-1}{k+1}}, }\quad
h_3(\tau) = \dstyle{ \left(  \frac{1}{\sqrt{2}} \frac{\eta(\tau)^2}{\eta(2\tau)} \right)^{\frac{k-1}{k+1}}. }
\end{align}
\end{itemize}
\end{thm}
\begin{cor}
\begin{align}
& S^J(1,k) =
\left( e^{ \frac{\pi i}{2} } \tau^{-1} \right)^{ (k-1)\frac{K}{2\kappa} } \frac{1}{\sqrt{2}}
\begin{pmatrix}
1  &  1   \\
1  &  -1  \\
\end{pmatrix},
\\[2mm]
& S^J(2,k) =
\left( e^{ \frac{\pi i}{2} } \tau^{-1} \right)^{ (k-1)\frac{K}{2\kappa} } \frac{1}{2}
\begin{pmatrix}
   1                       &             \sqrt{2}^{ 1+\frac{k-1}{k+1} }              &       1          \\
\sqrt{2}^{ 1-\frac{k-1}{k+1} }    &    0    &   - \sqrt{2}^{ 1-\frac{k-1}{k+1} }  \\
   1                       &           - \sqrt{2}^{ 1+\frac{k-1}{k+1} }              &       1
\end{pmatrix}.
\end{align}
\end{cor}
The matrix $S^J(1,k)$ can be easily calculated from theorems is Section \ref{sec:modular}.
The point is the expression of $S^J(2,k)$ which reveals its nice $k$ dependence. We tried to find this kind of simple $k$ dependence for $K=3,4$ by numerical calculation but failed.

\subsection{$S(K,k)$}
\label{subsec:S}
Define $\bar{J}_{\lambda, K}(z,u,\tau)$ by
\begin{align*}
\bar{J}_{\lambda, K}(z,u,\tau) &= e^{2\pi i(\kappa-2k)u} \,
\frac{ j_{-l}(z|\tau) }{ \hat{\delta}^{'k-1}},
\\
J_{\lambda, K} &= g_{l,K,k}\, \bar{J}_{\lambda, K},
\end{align*}
then
\begin{align*}
\bar{J}_{\mu, K} \Bigl( \frac z\tau ,u-\frac{z^2}{\tau}, -\frac 1\tau \Bigr)
 & = (-1)^{k-1}
   \, e^{ \frac{\pi i}{2} (k-1)\frac{K}{2\kappa} }
   \tau^{ -(k-1)\frac{K}{2\kappa} }
   \frac{2}{\sqrt{2\kappa}}
   \sum_{\lambda\in P^+_K} S(K,k)_{m,l}
   \: \bar{J}_{\lambda, K} (z ,u, \tau)
\\
 &= \sum_{\lambda\in P^+_K} S^{\bar{J}} (K,k)_{m,l} \:\bar{J}_{\lambda, K} (z ,u, \tau).
\end{align*}
The $S(K,1)=S^{\bar{J}}(K,1)=S^{J}(K,1)$ is the matrix of S-transformation for the characters of affine Lie algebra and has nice symmetry \cite{Ka}\cite{W}. Below we show the expressions of $S(K,k)$ for some lower $K$ which reflect this symmetry:
\begin{align*}
S(0,k) &= (-1)^{k-1} \sqrt{ \frac{\kappa}{2} },
\\[2mm]
S(1,k) &=
(-1)^{k-1} \sqrt{ \frac{\kappa}{2^2} }
\begin{pmatrix}
1  &  1   \\
1  &  -1  \\
\end{pmatrix},
\\[2mm]
S(2,k) &=
(-1)^{k-1} \sqrt{ \frac{\kappa}{2^3} }
\begin{pmatrix}
   1            &          1/s(k)         &       1               \\
   2 s(k)    &            0                &   - 2 s(k)        \\
   1            &        - 1/s(k)        &       1
\end{pmatrix},
\\[2mm]
S(3,k) &=
(-1)^{k-1} \frac{1}{2}
\sqrt{ \frac{\kappa}{1+ b' c'} }
\begin{pmatrix}
    1     &      c'      &      c'      &       1         \\
    b'    &      1       &     -1     &      -b'        \\
    b'    &     -1      &     -1     &        b'        \\
   1      &     -c'     &      c'      &      -1
\end{pmatrix}
\\[2mm]
&=
(-1)^{k-1}
\sqrt{ \frac{\kappa}{ 4(a^{\prime\prime 2}+ b'' c'')} }
\begin{pmatrix}
    a''     &      c''      &      c''      &       a''         \\
    b''    &      a''       &     -a''     &      -b''        \\
    b''    &     -a''      &     -a''     &        b''        \\
   a''      &     -c''     &      c''      &      -a''
\end{pmatrix},
\\
b' & =\frac{s(3)}{s(1)},\quad c'=\frac{s(k+1)}{s(k)},
\\
a'' & = s(1),\quad b''  =s(2k),\quad c''=s(k+1)\frac{s(1)}{s(k)},
\displaybreak[3]
\\[2mm]
S(4,k) &=
(-1)^{k-1}
\sqrt{ \frac{\kappa}{8 e^2} }
\begin{pmatrix}
    1     &      d      &      g       &        d        &       1         \\
    b     &      e       &     0       &      -e        &      -b         \\
    c    &       0      &      -2        &        0        &       c         \\
    b     &     -e       &     0       &       e        &      -b         \\
    1     &     -d      &      g       &      -d        &       1
\end{pmatrix},
\\[2mm]
b &= \frac{s(4)}{s(1)}, \quad d = \frac{s(k+1)}{s(k)}, \quad e =\frac{s(2)}{s(1)}, \quad g = \frac{1}{s(k)},
\\
c &= \frac{ s(k+1)s(k) }{ s(1)s(2) } P^{(k)}_2(k),
\\[2mm]
& \!\! P^{(k)}_2 (k) = m_{(2)} + \frac{ (1+q^2)(1-q^{2k}) }{ (1-q^{2(k+1)}) } m_{(1,1)},
\\[2mm]
&\qquad\quad = q^{2k} + q^{-2k} + \frac{ (1+q^2)(1-q^{2k}) }{
(1-q^{2(k+1)}) }.
\end{align*}
For $K=4$, we also have
$gc + 2 = 2 bd = 2 e^2$
which implies
$c = 2 s(k) \left( \left( \frac{s(2)}{s(1)} \right)^2 -1 \right)$.

In the rest of this subsection we demonstrate a way to calculate these $S(K,k)$ for $K=2$ case.

\noi {\it Step 1\ }: From Proposition \ref{prop:S-1} we can identify matrix elements to some extent:
\begin{align*}
S(2,k) =
\begin{pmatrix}
   a            &          c         &       a               \\
   b            &          0         &     - b             \\
   a            &         -c        &       a
\end{pmatrix}.
\end{align*}

\noi {\it Step 2\ }: From Proposition \ref{prop:S-2}, we have
\begin{align*}
S(2,k) &= (-2)^{k-1} \prod_{j=1}^{k} s(j)\,
\begin{pmatrix}
   1       &                                          &                                                             \\
            &      \frac{s(k)}{s(1)}       &                                                             \\
            &                                          &       \frac{ s(k)s(k+1) }{ s(1)s(2) }
\end{pmatrix}
\\[2mm]
& \qquad\qquad \times
\begin{pmatrix}
            1               &                1                    &                  1                     \\
  P^{(k)}_1(k)     &      P^{(k)}_1(k+1)       &        P^{(k)}_1(k+2)       \\
  P^{(k)}_2(k)     &      P^{(k)}_2(k+1)       &        P^{(k)}_2(k+2)
\end{pmatrix}
\begin{pmatrix}
   1       &                                              &                                                             \\
            &      \frac{s(k+1)}{s(k)}       &                                                             \\
            &                                               &          1
\end{pmatrix}
\\[2mm]
& = (-2)^{k-1} \prod_{j=1}^{k} s(j) \times
\begin{pmatrix}
                1                       &      \frac{1}{s(k)}          &                1                      \\
    \frac{s(2)}{s(1)}        &                  0                    &      -\frac{s(2)}{s(1)}    \\
                1                       &      -\frac{1}{s(k)}        &                1
\end{pmatrix}
\end{align*}
where we used
\[
 P^{(k)}_1(l) = 2\cos \left( \frac{l}{\kappa}\pi \right).
\]

\noi {\it Step 3\ }: We consider the constraint
\begin{align}
(S^{\bar{J}})^2 \times \left( e^{ -\frac{\pi i}{2} } \tau \right)^{ -(k-1)\frac{K}{2\kappa} }
=
(\mbox{phase factor}) \times \mbox{id}
\label{eqn:S^2=1}
\end{align}
which originates from the relation $S^4=1$ in $\Gamma$. Since we are dealing with a projective representation of $\Gamma$, there can be a phase factor.

From
\begin{align*}
\begin{pmatrix}
                1                       &      \frac{1}{s(k)}          &                1                      \\
    \frac{s(2)}{s(1)}        &                  0                    &      -\frac{s(2)}{s(1)}    \\
                1                       &      -\frac{1}{s(k)}        &                1
\end{pmatrix}^2
\propto \mbox{id},
\end{align*}
we have $\frac{1}{s(k)} \frac{s(2)}{s(1)}  = 2$. This with \eqref{eqn:S^2=1} implies
\begin{align*}
\frac{2}{\kappa} \left( (-2)^{k-1} \prod_{j=1}^{k} s(j) \right)^2 \times 4 = \pm 1.
\end{align*}
But clearly the R.H.S is $+1$ and
\begin{align*}
2^{k-1} \prod_{j=1}^{k} s(j) = \sqrt{ \frac{\kappa}{8} }.
\end{align*}

\subsection{$J_{\lambda, K}$ from $\bar{J}_{\lambda, K}$}
Here we show how to obtain $J_{\lambda, K}$, taking $K=2$ as an example.

\noi {\it Step 1\ }: Find $S^{\bar{J}} (\kappa,k)$: See Subsection \ref{subsec:S}.

\noi {\it Step 2\ }: Guess the transition matrix from its modular transformation:

Set $\bar{a}_{i,j}(\tau)$ by
\begin{align*}
\begin{pmatrix}
 \chi_{2\Lambda_0} \\
 \chi_{\Lambda_0 + \Lambda_1} \\
 \chi_{2\Lambda_1}
\end{pmatrix}
=
\begin{pmatrix}
   \bar{a}_{0,0}       &                  0              &       \bar{a}_{0,2}         \\
   0                  &           \bar{a}_{1,1}          &             0              \\
   \bar{a}_{2,0}       &                  0              &       \bar{a}_{2,2}
\end{pmatrix}
\begin{pmatrix}
 \bar{J}_{0.K} \\
 \bar{J}_{\bar{\Lambda}_1.K} \\
 \bar{J}_{2\bar{\Lambda}_1.K}
\end{pmatrix}.
\end{align*}
From modular transformation of $\chi_\lambda, \bar{J}_{\lambda,K}$, we have
\begin{gather*}
  \bar{a}_{0,0} = \bar{a}_{2,2}, \quad \bar{a}_{0,2} = \bar{a}_{2,0},
\\[5mm]
  \bar{a}_{0,0}(\tau+1) = e^{2\pi i \frac{k-1}{16(k+1)} } \,\bar{a}_{0,0}(\tau),
  \\
  \bar{a}_{2,0}(\tau+1) = e^{2\pi i \frac{9k+7}{16(k+1)} } \,\bar{a}_{2,0}(\tau),
  \\
  \bar{a}_{1,1}(\tau+1) = \bar{a}_{1,1}(\tau),
\\[5mm]
  \begin{pmatrix}
   \bar{a}_{0,0}(-1/\tau) \\
   \bar{a}_{2,0}(-1/\tau) \\
   \bar{a}_{1,1}(-1/\tau)
  \end{pmatrix}
  = \frac{1}{4} \left( e^{\frac{\pi i}{2}} \tau^{-1} \right)^{ -\frac{k-1}{2(k+1)} }
  \begin{pmatrix}
          2                  &                2                 &           2\sqrt{2} s(k)        \\
          2                  &                2                 &          - 2\sqrt{2} s(k)       \\
   2\sqrt{2} /s(k)     &    - 2\sqrt{2} /s(k)       &                0
  \end{pmatrix}
  \begin{pmatrix}
   \bar{a}_{0,0}(\tau) \\
   \bar{a}_{2,0}(\tau) \\
   \bar{a}_{1,1}(\tau)
  \end{pmatrix}.
\end{gather*}
To find $\bar{a}_{i,j}(\tau)$, introduce $h_i(\tau)$ by the base change
\begin{align*}
  \begin{pmatrix}
   h_{1}(\tau) \\
   h_{2}(\tau) \\
   h_{3}(\tau)
  \end{pmatrix}
  =
  \begin{pmatrix}
          1     &     1       &       0        \\
          1     &    -1       &       0       \\
          0     &     0       &   \sqrt{2} s(k)
  \end{pmatrix}
  \begin{pmatrix}
   \bar{a}_{0,0}(\tau) \\
   \bar{a}_{2,0}(\tau) \\
   \bar{a}_{1,1}(\tau)
  \end{pmatrix},
\end{align*}
then
\begin{align*}
  \begin{pmatrix}
   h_{1}(-1/\tau) \\
   h_{2}(-1/\tau) \\
   h_{3}(-1/\tau)
  \end{pmatrix}
 & = \left( e^{\frac{\pi i}{2}} \tau^{-1} \right)^{ -\frac{k-1}{2(k+1)} }
  \begin{pmatrix}
          1     &     0       &       0        \\
          0     &     0       &       1       \\
          0     &     1       &       0
  \end{pmatrix}
  \begin{pmatrix}
   h_{1}(\tau) \\
   h_{2}(\tau) \\
   h_{3}(\tau)
  \end{pmatrix},
\\[2mm]
  \begin{pmatrix}
   h_{1}(\tau+1) \\
   h_{2}(\tau+1) \\
   h_{3}(\tau+1)
  \end{pmatrix}
 & =
  \begin{pmatrix}
          0             &     e^{2\pi i \frac{k-1}{16(k+1)} }       &       0        \\
   e^{2\pi i \frac{k-1}{16(k+1)} }     &             0              &       0       \\
          0                       &                            0                           &       1
  \end{pmatrix}
  \begin{pmatrix}
   h_{1}(\tau) \\
   h_{2}(\tau) \\
   h_{3}(\tau)
  \end{pmatrix}.
\end{align*}
This transformation formulas looks similar to those of the Virasoro character appear in $c=1/2$ conformal field theory.
From these we can guess the form of $h_i$ as given in Theorem \ref{thm:K=1,2}.

\noi {\it Step 3\ }: Construct $\bar{J}_{\lambda,K}$ with the hypothetical transition matrix obtained in the last step:
\begin{align*}
\bar{J}_{0, 2} &= \frac{ \bar{a}_{0,0}\chi_{2\Lambda_0} - \bar{a}_{2,0}\chi_{2\Lambda_1} }{ \bar{a}_{0,0}^2 - \bar{a}_{2,0}^2 }
\\[2mm]
&= \frac{ \chi_{2\Lambda_0} + \chi_{2\Lambda_1} }{ 2h_1 }
      + \frac{ \chi_{2\Lambda_0} - \chi_{2\Lambda_1} }{ 2h_2 },
\\[2mm]
\bar{J}_{2\bar{\Lambda}_1, 2}
  &= \frac{ -\bar{a}_{2,0}\chi_{2\Lambda_0} - \bar{a}_{0,0}\chi_{2\Lambda_1} }{ \bar{a}_{0,0}^2 - \bar{a}_{2,0}^2 }
\\[2mm]
&= \frac{ \chi_{2\Lambda_0} + \chi_{2\Lambda_1} }{ 2h_1 }
      - \frac{ \chi_{2\Lambda_0} - \chi_{2\Lambda_1} }{ 2h_2 },
\\[2mm]
\bar{J}_{\bar{\Lambda}_1, 2} &= \frac{\chi_{\Lambda_0+\Lambda_1}}{ \bar{a}_{1,1} }
\\[2mm]
&= \frac{ \sqrt{2}s(k) }{h_3(\tau)} \chi_{\Lambda_0+\Lambda_1}.
\end{align*}
After normalization we obtain candidates for $J_{\lambda,K}$:
\begin{align*}
J_{0, 2} &= \bar{J}_{0, 2},
\\
J_{2\bar{\Lambda}_1, 2} &= \bar{J}_{2\bar{\Lambda}_1, 2},
\\
J_{\bar{\Lambda}_1, 2} &=
  \left( \frac{1}{\sqrt{2}} \right)^{\frac{k-1}{k+1}} \frac{ \chi_{\Lambda_0+\Lambda_1} }{ h_3 }.
\end{align*}

\noi {\it Step 4\ }: Check the defining relations for the $J_{\lambda,K}$'s:

Define a transition matrix $b_{\lambda, \mu} (\tau)$ by
\begin{align*}
J_{\lambda, K} (z,u,\tau)
= \sum_{\mu\in P^+_K} b_{\lambda, \mu} (\tau) \,\chi_{\mu, K} (z,u,\tau).
\end{align*}
We rewrite the defining condition (2) of Definition \ref{def:affineJack} as a first order differential equation for $b_{\lambda, \mu} (\tau)$
\begin{gather}
\sum_{\mu\in P^+_K} \,(k-1) \,b_{\lambda, \mu} \,\hat{\Delta} \chi_{\mu, K}
+ 2 (K + kh^\vee) \,(p\partial_p b_{\lambda, \mu}) \,\chi_{\mu, K} = 0,
\label{eqn:M'J}
\\
\hat{\Delta} = -\frac{1}{8\pi^2} \partial_{zz} - \frac{1}{2\pi^2}\partial_u\partial_\tau.
\end{gather}
Using the following lemma, we can check that this equation is satisfied by $b_{\lambda, \mu}$ given in Theorem \ref{thm:K=1,2}.
\begin{lem}
\label{lem:two}
\begin{align}
& \chi_{\Lambda_0 + \Lambda_1}(z,u,\tau) =
e^{ 4\pi iu } \vartheta_{2} (2z|\tau) \frac{ \eta(2\tau) }{ \eta(\tau)^2 },
\label{eqn:chi_1}
\\
& \frac{\hat{\Delta} \chi_{\Lambda_0 + \Lambda_1} }{ \chi_{\Lambda_0 + \Lambda_1} }
= 4p\partial_p \ln \frac{ \eta(\tau)^2 }{ \eta(2\tau) },
\\
& \chi_{2\Lambda_0}(z,u,\tau) - \chi_{2\Lambda_1}(z,u,\tau) =
e^{ 4\pi iu } \vartheta_{0} (2z|\tau) \frac{ \eta(\tau/2) }{ \eta(\tau)^2 },
\label{eqn:chi_2}
\\
& \frac{\hat{\Delta} ( \chi_{2\Lambda_0} - \chi_{2\Lambda_1} )}
{ \chi_{2\Lambda_0} - \chi_{2\Lambda_1} }
= 4p\partial_p \ln \frac{ \eta(\tau)^2 }{ \eta(\tau/2) },
\\
& \chi_{2\Lambda_0}(z,u,\tau) + \chi_{2\Lambda_1}(z,u,\tau) =
e^{ 4\pi iu } \vartheta_{3} (2z|\tau) \frac{ \eta(\tau) }{ \eta(\tau/2)\eta(2\tau) },
\label{eqn:chi_3}
\\
& \frac{\hat{\Delta} ( \chi_{2\Lambda_0} + \chi_{2\Lambda_1} )}
{ \chi_{2\Lambda_0} + \chi_{2\Lambda_1} }
= 4p\partial_p \ln \frac{ \eta(\tau/2)\eta(2\tau) }{ \eta(\tau) }.
\end{align}
\end{lem}
\begin{proof}
Let
\begin{align}
\Theta_{n,m} (z,u,\tau) &= e^{2\pi imu}\sum_{l \in \mathbb{Z} + n/2m } e^{2\pi i m(l^2\tau + lz)}
\\
&= e^{2\pi imu} \vartheta_{n,m} (z | \tau).
\end{align}
See Appendix \ref{appendix} for $\vartheta_i (v | \tau) \; (i=0,1,2,3)$.

For \eqref{eqn:chi_1}
\begin{align*}
& \vartheta_{2,4}(z|\tau) - \vartheta_{-2,4}(z|\tau) = i \vartheta_{1}(2z|2\tau),
\\
& \vartheta_{1}(2z|2\tau) =
\vartheta_{1}(z|\tau) \vartheta_{2}(z|\tau) \frac{ \eta(2\tau) }{ \eta(\tau)^2 }.
\end{align*}
For \eqref{eqn:chi_2}
\begin{align*}
\Bigl(
\Theta_{1,4} - \Theta_{-1,4} - \left( \Theta_{3,4} - \Theta_{-3,4} \right)
\Bigr) (z,u,\tau)
&= i e^{8\pi iu} \vartheta_1 (z|\tau/2)
\\
&= i e^{8\pi iu} \vartheta_1 (z|\tau) \vartheta_0 (z|\tau)
\frac{ \eta(\tau/2) }{ \eta(\tau)^2 }.
\end{align*}
For \eqref{eqn:chi_3}
\begin{align*}
& \Bigl(
\Theta_{1,4} - \Theta_{-1,4} + \Theta_{3,4} - \Theta_{-3,4}
\Bigr) (z,u,\tau+1)
\\
&= e^{\pi i/8}
\Bigl(
\Theta_{1,4} - \Theta_{-1,4} - \Theta_{3,4} + \Theta_{-3,4}
\Bigr) (z,u,\tau)
\\
&= e^{\pi i/8} i e^{8\pi iu} \vartheta_1 (z|\tau/2).
\end{align*}
\end{proof}

The defining condition (1) of Definition \ref{def:affineJack} can be easily checked.
This completes the procedure to find $J_{\lambda, K}$ for $K=2$.

\noi
{\it Remark}\;
For the case of level one, the equation \eqref{eqn:M'J} reduces to
\begin{align*}
\frac{\hat{\Delta} \chi_{\Lambda_i} }{ \chi_{\Lambda_i} } &= 2p\partial_p \ln \eta(\tau).
\end{align*}

\noi
{\it Remark}\:
It is possible to find $b_{\lambda, \mu}$ by solving \eqref{eqn:M'J} for level one and two.
In fact this is done in \cite{FSV2} and some more result is obtained.

\bigskip
As byproducts, we have this lemma.
\begin{lem}
\label{lem:g}
\begin{itemize}
\item[]
\item[(1)]\quad
$\frac{{\displaystyle g_{k+1,2,k} }}{{\displaystyle g_{k,2,k} }} \times \sqrt{2}^{ 1+\frac{k-1}{k+1} } s(k) = 1,$
\item[(2)]\quad
$g_{m,K,k} = g_{p-m,K,k}.$
\end{itemize}
\end{lem}
\noi
Direct proof is also available from \eqref{eqn:g}.

\appendix
\section{Theta functions}
\label{appendix}
Here we summarize the definition of theta functions $\vartheta_i (v | \tau) \; (i=0,1,2,3)$.
\begin{align*}
\vartheta_1 (v | \tau) &= i\sum_{n\in\Z} (-1)^n e^{i\pi\tau (n-\frac 12)^2} e^{i\pi v(2n-1)},
\\
\vartheta_2 (v | \tau) &= \sum_{n\in\Z} e^{i\pi\tau (n-\frac 12)^2} e^{i\pi v(2n-1)},
\\
\vartheta_3 (v | \tau) &= \sum_{n\in\Z} e^{i\pi\tau n^2} e^{i\pi v \,2n},
\\
\vartheta_0 (v | \tau) &= \sum_{n\in\Z} (-1)^n e^{i\pi\tau n^2} e^{i\pi v \,2n}.
\end{align*}

\vspace{1cm}
\noi
{\Large\it Acknowledgements}\\
The author would like to thank M. Jimbo, A. Kirillov Jr., K. Mimachi, M. Nishizawa and J. Shiraishi for fruitful discussions and helpful comments.
He is a Research Fellow of the Japan Society of the Promotion of Science. This work is supported by the Grant-in-Aid for Scientific Research from Ministry of Education, Culture, Sports, Science and Technology of Japan No.13-04575.


\newpage


\begin{thebibliography}{Lampoxlong}

\bibitem[BeFe]{BeFe}
D. Bernard and G. Felder,
Fock Representations and BRST cohomology in $SL(2)$ Current Algebra,
{\sl Comm. Math. Phys.} {\bf 127}, 145-168 (1990).

\bibitem[C]{C}
I. Cherednik,
Double affine Hecke algebras and Macdonald conjectures,
{\sl Annals of Math.} {\bf 141} (1995), 191-216.

\bibitem[CMS]{CMS}
J. van Diejen and L. Vinet eds.,
Calogero-Moser-Sutherland Models,
Springer-Verlag, New York, (2000).

\bibitem[EK]{EK}
P. Etingof and A. Kirillov Jr.,
Macdonald's polynomials and representations of quantum groups,
{\sl Math. Res. Lett.} {\bf 1} (1994), no. 3, 279--296;
Representation-theoretic proof of the inner product and symmetry identities for Macdonald's polynomials,
{\sl Compositio Math.} {\bf 102} (1996), no. 2, 179--202.;
On the affine analogue of Jack's and Macdonald's polynomials,
{\sl Duke Math. J.} {\bf 78} (1995), no. 2, 229--256.

\bibitem[FF]{FeFr}
B. Feigin and E. Frenkel,
{\sl Representations of Affine Kac-Moody Algebras and Bosonization} in
{\sl Physics and Mathematics of Strings},
World Scientific, Singapore (1990).

\bibitem[Fe]{Fe}
G. Felder,
BRST Approach to Minimal Models,
{\sl Nucl. Phys.} {\bf B317} (1989) 215-236,
Erratum {\it ibid}. {\bf B324}, 548 (1989).

\bibitem[FeSi]{FeSi}
G. Felder and R. Silvotti,
Modular covariance of minimal model correlation functions,
{\sl Comm. Math. Phys.} {\bf 123} (1989), no. 1, 1--15.

\bibitem[FV]{FV}
G. Felder and A. Varchenko,
Integral representations of solutions of the elliptic Knizhnik-Zamolodchikov-Bernard equations,
{\sl Internat. Math. Res. Notices} (1995) no. 5, 221-233.

\bibitem[FSV1]{FSV}
G. Felder, L. Stevens and A. Varchenko,
Modular transformations of the elliptic hypergeometric functions, Macdonald polynomials, and the shift operator,
math.QA/0203049.

\bibitem[FSV2]{FSV2}
G. Felder, L. Stevens and A. Varchenko
Elliptic Selberg Integrals and Conformal Blocks,
math.QA/0210040.

\bibitem[Ka]{Ka}
V. Kac
{\sl Infinite Dimensional Lie Algebras}, 3rd ed.,
Cambridge Univ. Press (1990).

\bibitem[Ki]{Ki}
A. Kirillov Jr.,
On an inner product in modular tensor categories,
{\sl J. Amer. Math. Soc.} {\bf 9} (1996), no. 4, 1135;
On inner product in modular tensor categories. II. Inner product on conformal blocks and affine inner product identities,
{\sl Adv. Theor. Math. Phys.} {\bf 2} (1998), no. 1, 155--180.

\bibitem[M1]{M1}
I. Macdonald,
{\sl Symmetric Functions and Hall Polynomials}, 2nd ed.,
The Clarendon Press, Oxford University Press, New York (1995).

\bibitem[M2]{M}
I. Macdonald,
{\sl Symmetric Functions and Orthogonal Polynomials},
University Lecture Series, 12. (American Mathematical Society, Providence, RI, 2001).

\bibitem[N]{N}
M. Noumi,
Macdonald's symmetric polynomials as zonal spherical functions on some quantum homogeneous spaces,
{\sl Adv. Math.} {\bf 123} (1996), no. 1, 16--77.

\bibitem[W1]{W}
M. Wakimoto,
{\sl Infinite-dimensional Lie algebras}, Translations of Mathematical Monographs, 195.
(American Mathematical Society, Providence, RI, 2001).

\bibitem[W2]{W2}
M. Wakimoto,
Fock representations of affine lie algebra {$A^{(1)}_1$}.
{\sl Comm. Math. Phys.}, {\bf 104} (1986), 605--609.

\end{thebibliography}
\end{document}